\documentclass[12pt]{amsart}
\usepackage[latin1]{inputenc}
\usepackage{amsmath}
\usepackage{amsthm}
\usepackage{amsfonts}
\usepackage{amscd}
\usepackage{amssymb}
\usepackage{mathrsfs}
\usepackage{graphicx}
\usepackage{psfrag}
\usepackage{tikz}
\usepackage[all]{xy}
\usepackage{float}

\newcommand{\forget}[1]{{}}
\def\ca{{\mathcal A}}

\def\ch{{\mathcal H}}

\def\cl{{\mathcal L}}

\def\bn{{\mathbb N}}

\def\br{{\mathbb R}}

\def\bt{{\mathbb T}}
\def\bz{{\mathbb Z}}

\def\z{\zeta}
\def\x{\xi}
\def\r{\rho}

\def\f{\varphi}

\def\Tr{\mathrm{Tr}}

\newtheorem{Thm}{Theorem}[section]
\newtheorem{Cor}[Thm]{Corollary}
\newtheorem{Prop}[Thm]{Proposition}
\newtheorem{Lemma}[Thm]{Lemma}

\newtheorem{Conjecture}[Thm]{Conjecture}
\newtheorem{Dfn}[Thm]{Definition}

\newtheorem{exmp}[Thm]{Example}

\newtheorem{numbering}[Thm]{Numbering}
\newtheorem{Co}[Thm]{Comment}

\theoremstyle{remark}
\newtheorem{rem}[Thm]{Remark} 
 
\addtocounter{section}{-1}

\begin{document}

 \title[Spectral triples and the geometry of fractals]{Spectral triples and the geometry of fractals }

 \author{Erik Christensen, Cristina Ivan, Elmar Schrohe}

\address{  
 Department of Mathematics, University of Copenhagen,   \newline
\indent DK-2100 Copenhagen, Denmark}
 \email{echris@math.ku.dk}

\address
{ 
UT MD Anderson Cancer Center \newline
\indent Smith Research Building \newline
\indent Houston, Texas, 77054, USA}
\email{civan@mdanderson.org}

 \address{   
Institut f\"ur Analysis, Leibniz Universit\"at Hannover,  \newline
\indent   30167 Hannover, Germany }
 \email{schrohe@math.uni-hannover.de}

 \date{\today}



  \begin{abstract}
We construct spectral triples for the Sierpinski gasket as infinite sums of 
unbounded Fredholm modules associated with the holes in the gasket 
and investigate their properties. For each element in the K-homology group we find a representative induced by one of our spectral triples. Not all of these triples, however, will have the right geometric properties.  If we want the metric induced by the spectral triple to give the geodesic 
distance, then we will have to include  a certain minimal family of unbounded Fredholm modules.
If we want the eigenvalues of the associated generalized Dirac 
operator to have the right summability properties, then we get limitations
on the number of summands that can be included. 
If we want the Dixmier trace of the spectral triple to coincide with a multiple of the Hausdorff measure, then we must impose conditions on the distribution of the summands over the gasket.  
For the elements of a large subclass of the K-homology group, however, the representatives are induced by triples 
having the desired geometric properties. We finally show that the same techniques can be applied to the Sierpinski 
pyramid.  
\end{abstract}
\maketitle

\section{Introduction}
In his noncommutative geometry program Alain Connes employs ideas from operator algebras to analyze singular spaces for which the classical tools of geometric analysis fail.
One of the basic structures in this theory is that of a spectral triple 
$(\ca,\ch,D)$, consisting of an algebra $\ca$ of bounded operators on 
a Hilbert space $\ch$ and an unbounded selfadjoint operator $D$ on $\ch$.
In this picture, the space is replaced by the algebra $\ca$ -- in the simplest 
cases an algebra of sufficiently smooth functions on the space --
while the geometry is encoded in the operator $D$, which is required to have a  
compact resolvent and bounded commutators with the elements of $\ca$.
In \cite{Co2, Co4}   Connes proves with some relevant examples that his program may be used to study fractals. 
In \cite{La1, La2} Michel Lapidus investigates in many different ways the possibility of developing a noncommutative fractal geometry. 
In \cite{CI} Christensen and Ivan construct spectral triples
for approximately finite dimensional C*-algebras and then apply
this result to the special case of the continuous functions on the
Cantor set. In \cite{CIL} Christensen, Ivan and Lapidus construct
a spectral triple associated to the Sierpinski gasket. 
It encodes the geometry in that it recovers the geodesic distance, the Hausdorff dimension
and the Hausdorff measure; moreover, it gives a non-trivial
element in the K-homology group of the gasket. 

In this paper, our main interest is to determine which K-homology elements we can obtain by constructions of this type.
At the same time, we extend the analysis to the case of the Sierpinski pyramid. 
In both cases, a crucial role is played by the sets {\it HSG} and {\it HSP} of holes in the gasket and
non-horizontal holes in the pyramid, respectively. 
In fact, the $K_1$ groups are just the free abelian groups 
\begin{displaymath}
\underset{HSG}{\oplus}\bz \, \text{ and }\, \underset{HSP}{\oplus} \bz,
\end{displaymath}
and the K-homology groups $K^1$ are the dual groups with respect to
$\bz,$
\begin{displaymath}
\underset{HSG}{\Pi}\bz \, \text{ and }\, \underset{HSP}{\Pi}\bz.
\end{displaymath}
A K-homology element can therefore be identified with a sequence of integers indexed by {\it HSG} and {\it HSP},
respectively.

We show that it is possible to
obtain any element in the K-homology group from our construction, but
the only geometric structure which is preserved by all such spectral
triples is the geodesic distance. If we want to have a spectral triple of this type
which gives the right Hausdorff dimension, then there is some limitation 
on the growth of the associated sequence representing the K-homology element. 
Whenever this sequence is bounded, we find a spectral triple with the 
right metric and the right dimension. There is room
for some unbounded sequences too, but we are not able to say exactly
which ones we can obtain. 
The Hausdorff measure is dominated by the volume measure induced by any of our spectral triples. 
In general, however, the two measures are not proportional.  We find sufficient conditions for this to be true.  
These investigations show that inside our set-up there are bounds on
the numbers as well as on the distributions of summands we use in
forming our spectral triples. In the other direction we can prove that
the spectral triple, which we denote $ZGT, $ and which represents the
$0-$element in the K-homology group, actually has a minimality
property. We do this by showing that if just any summand -- in the sum
of unbounded Fredholm modules giving this spectral triple -- is left
out, then the geodesic distance will not be the metric coming from
this spectral triple.

\noindent{\it Acknowledgment.} We are grateful to Malek Joumaah for the drawings.

\section{Sierpinksi Gasket: Constructions, K-theory and K-homology}
There are two basic 
procedures, which both produce the Sierpinski gasket; one is based on 
continued cuttings and the other on continued extensions
of graphs. The first steps in the cutting procedure are shown in the
Figure \ref{SG1} below. 
\begin{figure}[H] 
\includegraphics[scale=0.2]{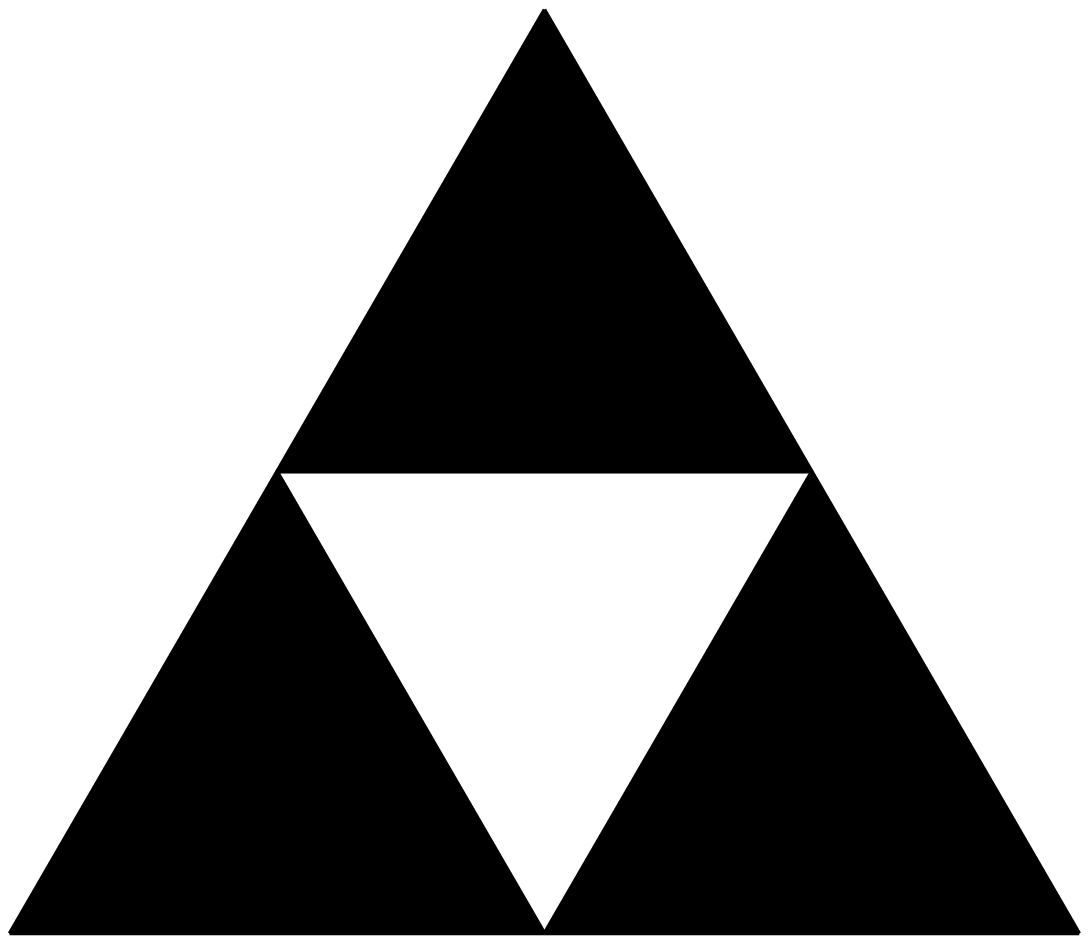}\includegraphics[scale=0.2]{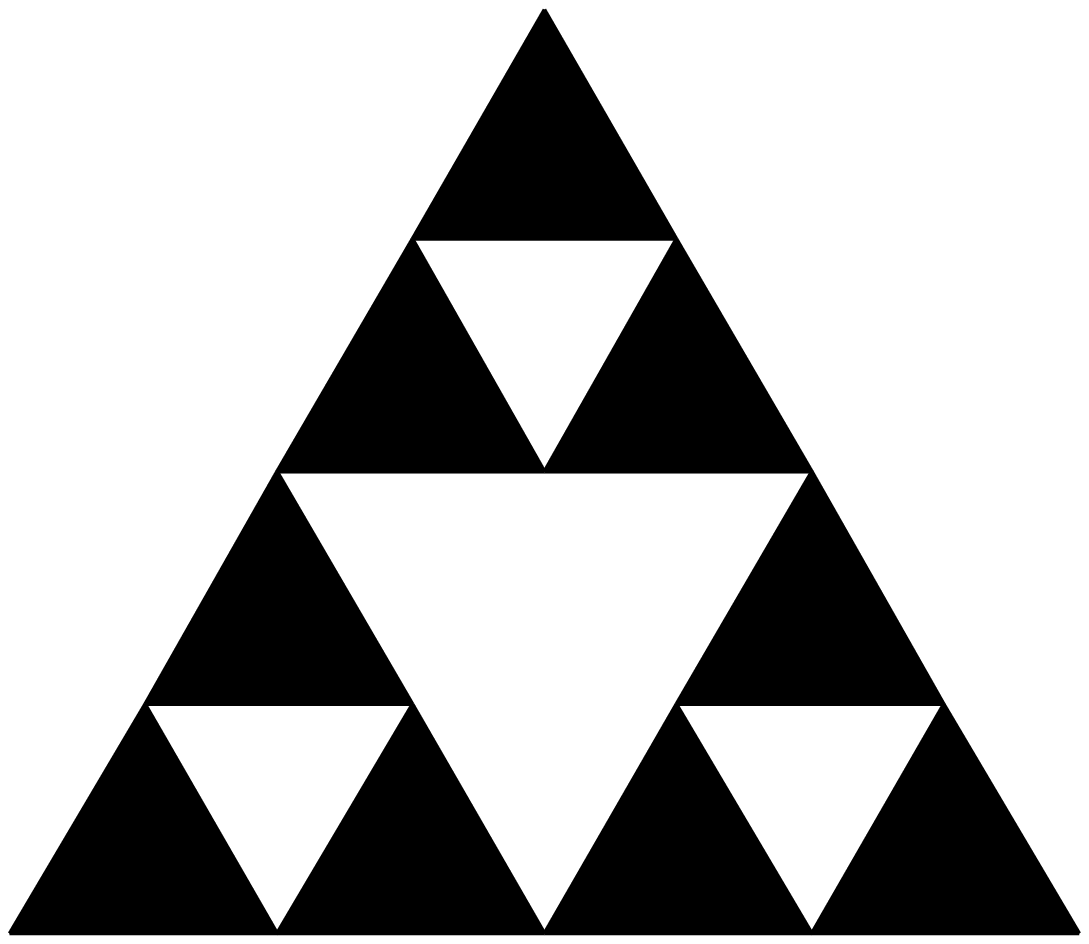}\includegraphics[scale=0.2]{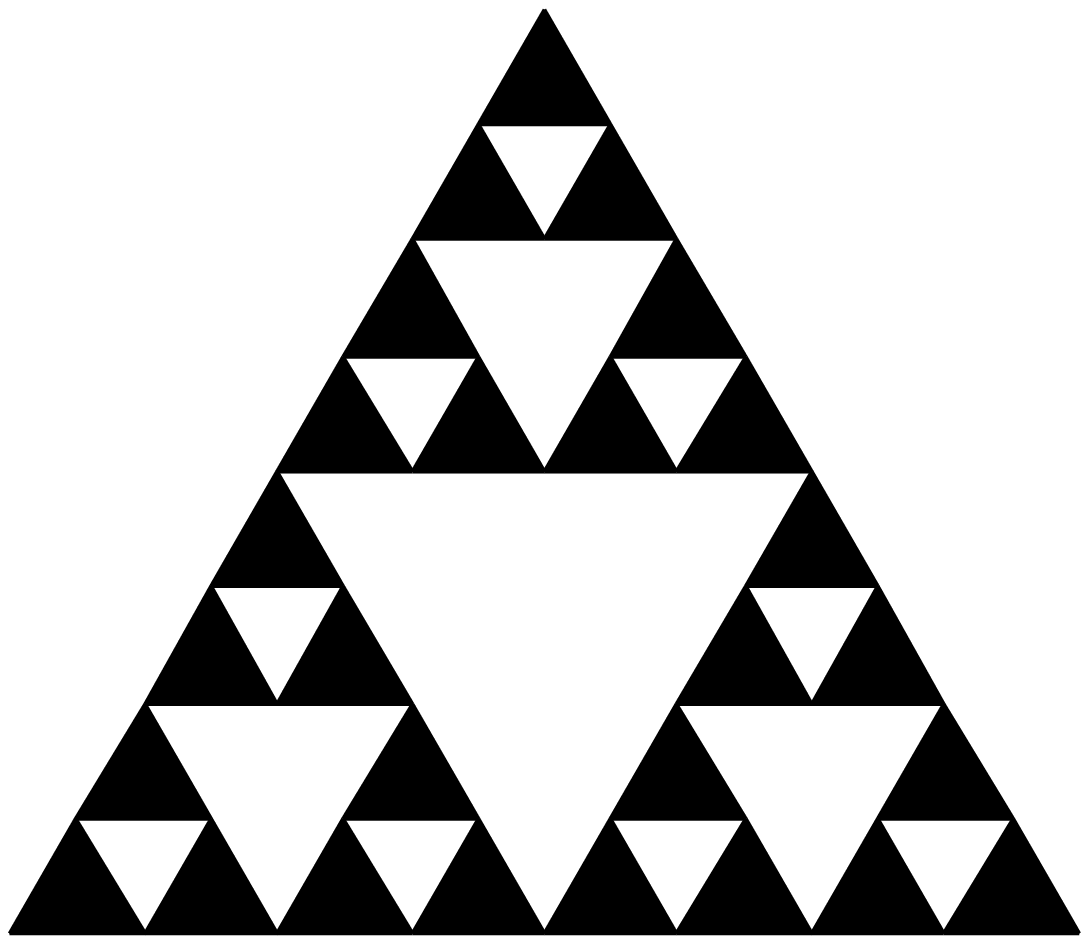}
\caption{\label{SG1}} \end{figure}  
\indent The Sierpinski Gasket is just the intersection of
all these sets, and a compactness argument shows that this is a non-empty
compact subset of the plane. 
The inductive construction procedure shows that there are many {\em holes}, i.e. bounded components of the
complement of the gasket, and once a hole has been added it will
remain undisturbed during the following steps of the construction. \\
\indent The first steps in the extension procedure are shown in the Figure \ref{SG2} below. 
\begin{figure}[H] \includegraphics[scale=0.2]{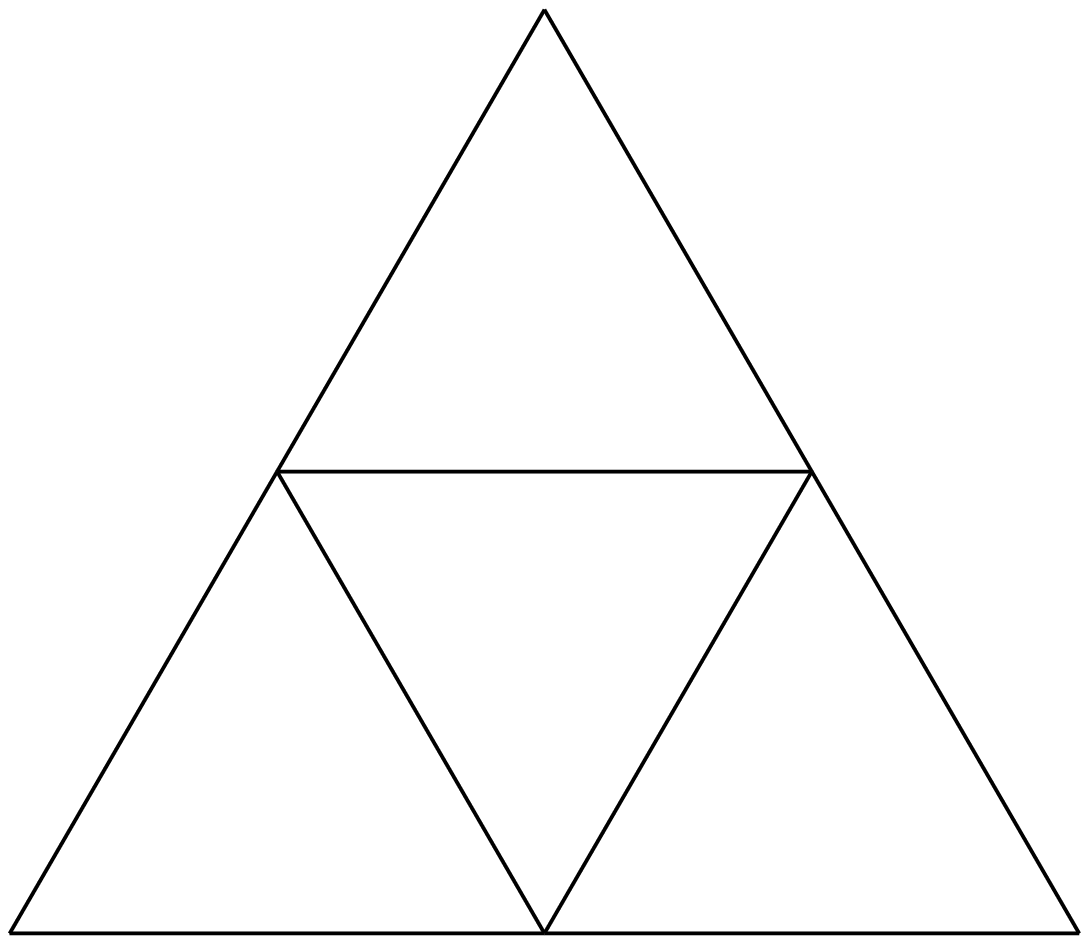}\includegraphics[scale=0.2]{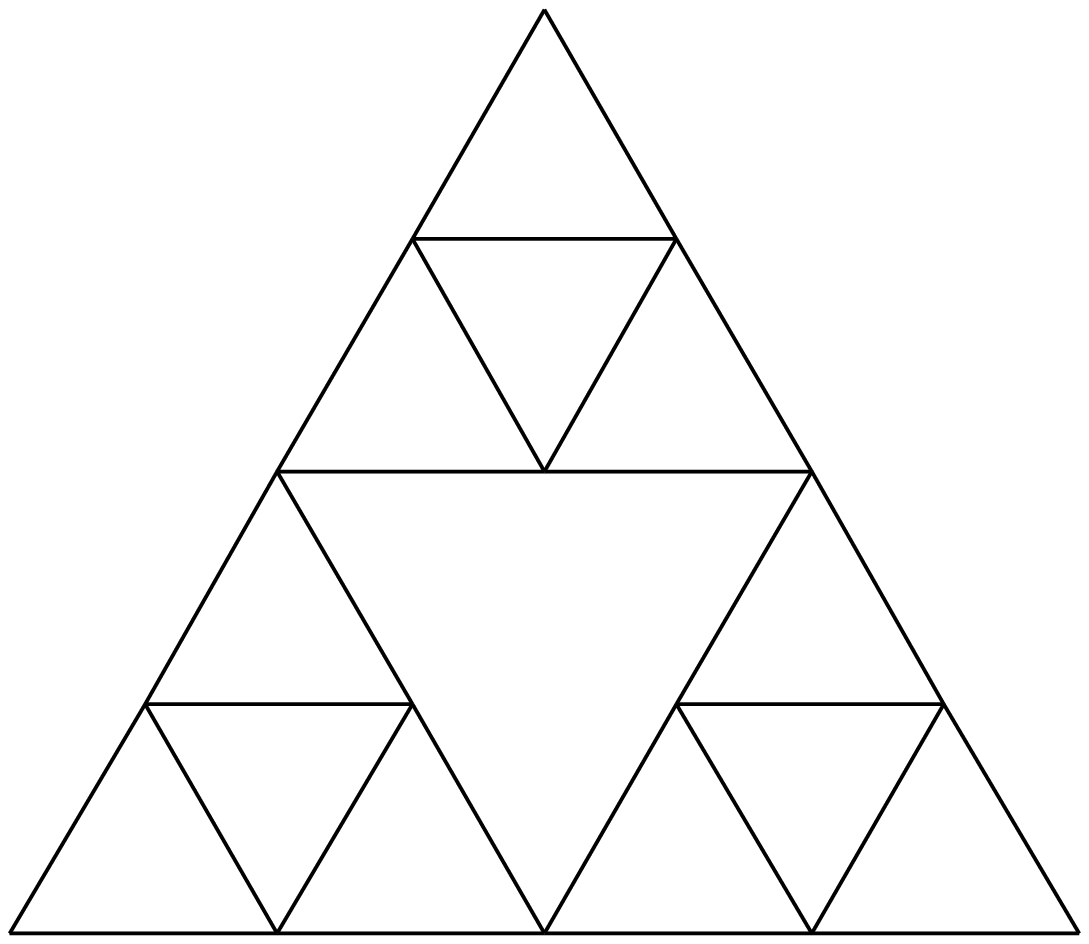} \includegraphics[scale=0.2]{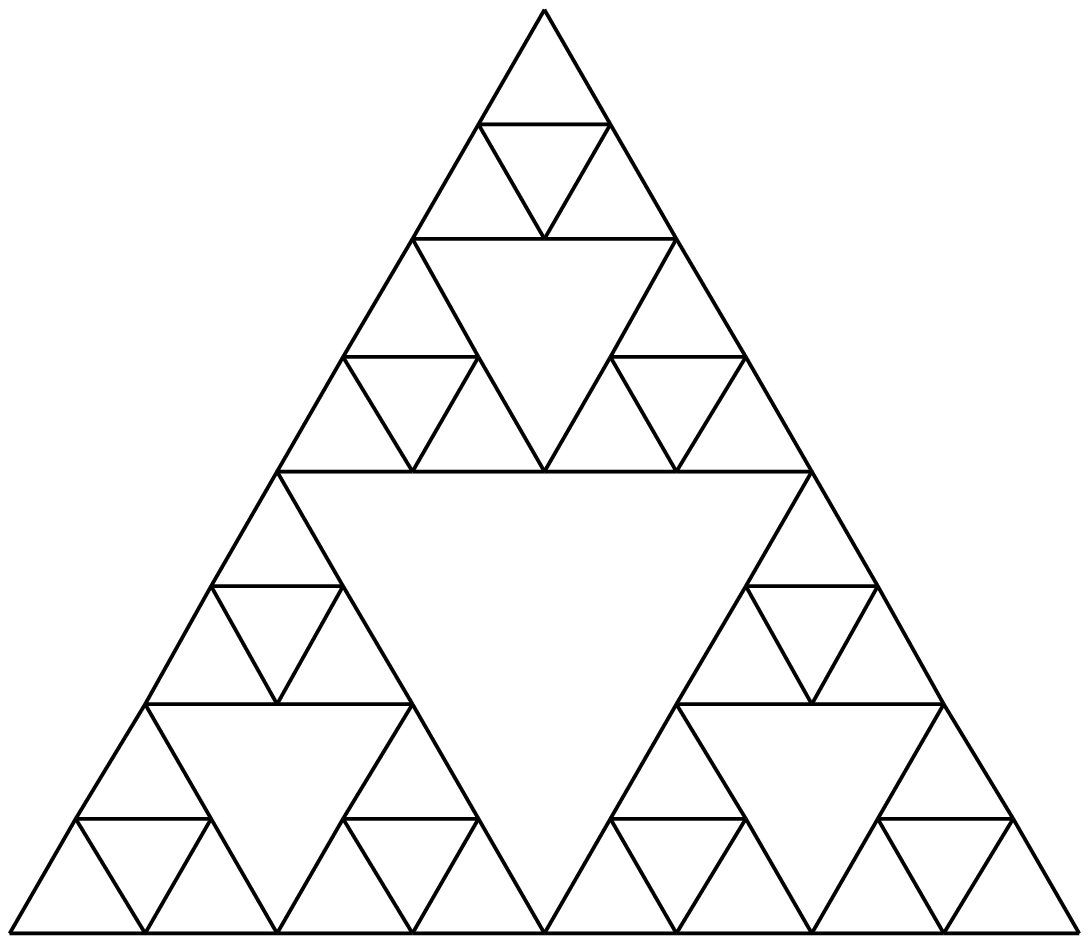}\caption{\label{SG2}}\end{figure}
\indent In this case the gasket is obtained as the closure of the union 
of all these sets, which may be thought of as planar graphs.\\
\indent 
The $K_1$ and the first Steenrod K-homology group $K^1$ of the gasket are computable from general results on the K-theory and K-homology of planar sets.  Let $X$ be a non-empty compact Hausdorff space and $\check{H}^1(X)$ the first {\v C}ech cohomology group of $X$, i.e. $\check{H}^1(X)$ is the quotient of the multiplicative group $C(X)^{-1}$ of continuous non-vanishing functions on $X$, by its component of unity. Define a character homomorphism $K_1(C(X)) \rightarrow \check{H}^1(X)$ by assigning to a unitary in $M_k(C(X))$ -- which is a continuous, unitary-valued function on $X$ -- the class of the function which is obtained as the pointwise determinant.  If $X$ is a non-empty, compact subset of the plane then this is an isomorphism (see for example \cite{HR}).\\
\indent The following description of $\check{H}^1(X)$ makes it clear that $K_1(C(X))$ is free abelian.
\begin{Thm}(\cite{HR})
Let $\{ \lambda_1, \lambda_2, \dots \}$ be a sequence of points in $\mathbb{C} \backslash X$ with precisely one 
$\lambda_j$ in each {\it hole}, i.e. bounded component of the complement of $X$, 
and none in the unbounded component of the complement. The group $\check{H}^1(X)$ is freely generated as an abelian group 
by the homotopy classes of the nowhere zero complex functions $z-\lambda_j$ on $X$.
\end{Thm}

For the first Steenrod K-homology group of $X$, $K^1(C(X))$, there are two equivalent descriptions.   The first is due to Brown, Douglas, Fillmore (see \cite{BDF}) who realize $K^1(C(X))$ as the abelian group of unitarily equivalent extensions of C*-algebras of the form 
\begin{displaymath}
0 \longrightarrow \mathcal{K} \longrightarrow A \longrightarrow C(X) \longrightarrow 0
\end{displaymath}
where $\mathcal{K}$ is the C*-algebra of compact operators on  a separable complex Hilbert space.  
The second description is due to  Kasparov \cite{Ka2} (see also \cite{Co1, Co2}) who realizes $K^1(C(X))$,    
as the abelian group of stable homotopy classes of bounded Fredholm modules over $C(X)$. \\

\indent If $X$ is a compact planar set 
\begin{displaymath}
K^1(C(X)) \cong \text{Hom}(K_1(C(X), \mathbb{Z}) \quad (\text{\cite{BDF}, \cite{Do}}). 
\end{displaymath}
In Kasparov's picture of $K^1$ the isomorphism between $K^1(C(X))$ and  $\text{Hom}(K_1(C(X), \mathbb{Z})$  is given by the following index map (\cite{Ka1},\cite{Co2}):

Let $(H (\pi), F)$  be an (odd bounded)  Fredholm module over $C(X)$ and let $P=(I+F)/2$.  Let $u$ be a unitary in $C(X)$.  Then the operator $P\pi(u)P$
from $PH$ to itself is a Fredholm operator.  An additive map from  $K_1(C(X))$ to $\mathbb{Z}$ is determined by 
\begin{displaymath}
\Phi_{(H,F)}([u])= - \text{Index }(P\pi(u)| PH) .
\end{displaymath}
The index map $\Phi_{(H,F)}$ only depends upon the class of $(H,F)$.

\indent Baaj and Julg have shown in \cite{BJ} that every class of bounded Fredholm modules over the algebra $C(X)$ contains a Fredholm module constructed from an unbounded one (see also \cite{Co2}, IV, Appendix A). Precisely, any unbounded Fredholm module $(H, D)$ associated to $C(X)$ defines a bounded Fredholm module $(H, F),$ where $F$ is the self-adjoint unitary coming
from the polar decomposition of the Dirac operator $D$.  The symmetry $F$ equals $2P-I$ where $P$ is the spectral projection $E([0,\infty[)$ for $D$.  Every bounded Fredholm module over $C(X)$ is operator homotopic to one obtained from the foregoing construction.\\
\begin{rem}
The above description of $K_1(C(X))$ for any compact planar set $X$ identifies $\text{Hom}(K_1(C(SG)), \mathbb{Z})$, and thus $K^1(C(SG))$, with the product group $\Pi \mathbb{Z}$, with one factor for each hole in $SG$. 
\end{rem}

\section{Some Notation and Conventions}  \label{OST}

The basis for our constructions is Section 8 of
\cite{CIL},  and we assume some familiarity with these results. There a spectral triple for the Sierpinski Gasket is
constructed as an example. We will call this triple {\em the old gasket triple}. % and use the acronym $OGT$ as its name in this article. 
%We will not repeat  the construction here, but 
Let us shortly recapture some of the definitions
 and names. 
%we used in the investigation of the $OGT$. 
The basis for the construction is an equilateral triangle of
 circumference $2\pi$ which is placed as in Figure \ref{SG2}
 such  that one vertex is pointing upwards. This triangle is denoted
 $\Delta_{0,1},$ and for any natural number $n$ the gasket contains
 $3^n$ triangles $\Delta_{n,j}, \, 1 \leq j \leq 3^n, $ of circumference
 $2\pi/(2^n)$ which are all scaled
 and translated copies of $\Delta_{0,1}$ as shown in Figure \ref{SG2}.
Each of the triangles $\Delta_{n,j}$ is then treated as a circle with radius $2^{-n}. $ In 
\cite{CIL}, Theorem 2.4, we investigated the standard spectral triple 
for a circle of radius $2^{-n},$ and we showed how this triple can be
transformed into an unbounded Fredholm module for $C(SG).$ We will
let $UFM(\Delta_{n,j})$ denote this unbounded Fredholm module.
The direct sum of these unbounded Fredholm modules 
over all the pairs $(n, j)$ then gives the old spectral triple. 
Later, the authors discussed
other possibilities for a spectral triple for the Sierpinski
Gasket. It is clear that the upside down triangles which form the
boundaries of the 
bounded components of the complement of the gasket 
must contain nearly
the same information as the triangles $\Delta_{n,j}.$ The triangle $\Delta_{0,1} $ then only will
appear indirectly as a boundary in a completion of the union of all
the upside-down triangles, but this is
not a serious problem, since we can just add the module corresponding
to  the outer triangle to the direct sum of the 
modules associated
to all the upside down triangles.
To be more precise we will also introduce a
numbering system for the upside down triangles. The central one is
denoted $\nabla_{1,1}; $ it is of circumference $\pi.$ 
 All the upside down triangles are then numbered by a pair $(m,k) $
 with $m, k$   natural numbers such that $1 \leq k  \leq 3^{m-1},$ 
  and they are denoted $\nabla_{m,k}.$ 
\indent With $HSG$ denoting the countable set that enumerates the holes in the Sierpinski Gasket, i.e  
\begin{displaymath} 
HSG \, := \, \left \{ (m,k)\,: \, m \in \bn , \, 1 \leq k \leq 3^{m-1} \,\right \}. 
\end{displaymath}
we may write 
\begin{align*}
K_1(C(SG))& = \oplus_{(m,k) \in HSG} \bz,\\
K^1(C(SG)) & =\mathrm{Hom}(K_1(C(SG)), \bz)=\Pi_{(m,k) \in HSG} \bz .
\end{align*}

\section{On a family of spectral triples representing any element in $K^1(C(SG))$ }
Just as for the triangles $\Delta_{n,j}$ we can
construct an unbounded Fredholm module, 
\begin{displaymath}
UFM(\nabla_{m,k})=(A, H_{\nabla_{m,k}}, D_{\nabla_{m,k}}),
\end{displaymath}
for the Sierpinski gasket by 
parameterizing  $\nabla_{m,k} $ on a circle of radius $2^{-m}$ following  \cite{CIL}, Definition 8.1. 
For each of the triangles $\Delta_{n,j}$, respectively $\nabla_{m,k}$,
we can also construct an unbounded Fredholm module by reversing the
orientation in the parameterization of the triangles. 
These modules will be denoted $\overline{UFM(\Delta_{n,j})}$ and $\overline{ UFM(\nabla_{m,k})},$ respectively.
\begin{Prop} \label{KUFM}
Let $u$ be a unitary in $C(SG),$  $[u]$ its class in \newline
 $K_1(C(SG))$.
For each of the triangles $\Delta_{n,j} ,$ respectively $\nabla_{m,k},$ denote the
winding number for the restriction of $u(z)$ to this triangle by $w_{\Delta_{n,j}}(u), $ respectively $w_{\nabla_{m,k}}(u).$ Then 
\begin{itemize} 
\item[(i)] $\Phi_{UFM(\Delta_{n,j})}([u]) = w_{\Delta_{n,j}}(u).$
\item[(ii)] $\Phi_{UFM(\nabla_{m,k})}([u]) = w_{\nabla_{m,k}}(u).$
\item[(iii)]$ \Phi_{UFM(\Delta_{n,j})}([u]) = 
\sum_{\nabla_{m,k} \subset \Delta_{n,j}}  w_{\nabla_{m,k}}(u).$
\item[(iv)] The  element  $ \Phi_{UFM(\Delta_{n,j})}$ 
in the group  \begin{displaymath} 
\mathrm{Hom}(K_1(SG), \bz) =
  \Pi_{(m,k) \in HSG}\bz
\end{displaymath} 
is the function $f_{\Delta_{n,j}} : HSG \to \bz$ 
given by 
\begin{displaymath}
f_{\Delta_{n,j}}(m,k) = \begin{cases} 1 \text{ if } \nabla_{m,k} \subset 
\Delta_{n,j} \\ 0 \text{ else }. \end{cases} \end{displaymath}
\end{itemize}   
\end{Prop}

\begin{proof}
(i) and (ii) are consequences of 
the result that for
the unit circle $\bt,$  the Hilbert space $H = L^2(\bt)$ and 
the projection $P_+$ onto $H_+ = \overline{\text{span}}\{z^n\, 
\big|\, n\geq 0 \},$  the winding number of $u$ 
on $\bt$ is the opposite integer to the index of $P_+M_u|H_+,$
 where $M_u$ is the multiplication operator induced by $u.$

For (iii) use the uniform continuity of $u$: From  a certain natural number $n_0$ on we have for any $n, m \geq n_0$
 that the winding 
numbers of $u$ around any $\Delta_{n,j}$ and any $\nabla_{m,k}$ 
vanish for $m,n\geq n_0.$ Let us next compute $ \Phi_{UFM(\Delta_{0,1})}([u]).$ 
By the  same argument as above, it equals the winding number of $u$ around 
$\Delta_{0,1}.$ On the other hand this is the sum of the winding
numbers over the four triangles $\Delta_{1,j}, \,\, 1 \leq j \leq 3,$ 
and $\nabla_{1,1}.$ For each of $\Delta_{1,j} $ 
we repeat the subdivision until we reach the level $n_0$ from where on all winding numbers of $u$
 vanish. We are left with the sum of all the winding
 numbers of $u$ over the triangles $\nabla_{m,k},\,1 \leq m
 \leq n_0, \, 1 \leq k \leq 3^{m-1},$ and obtain (iii) 
in the case where $n=0$ and $j=1.$ The general case then
follows by an analogous argument applied to the triangle
$\Delta_{n,j}. $
 
(iv) is just a reformulation of (iii) which is 
suitable for the  computations to come. 
\end{proof}

We recall from \cite{CIL} that we can perform  infinite-direct 
sums of unbounded Fredholm  modules of the types
\begin{displaymath}
UFM(\nabla_{m,k}) \quad \text{or} \quad \overline{UFM(\nabla_{m,k})}
\end{displaymath}
 and get an unbounded Fredholm
module as a result. 

We will now construct a spectral triple $ZGT$ for the Sierpinski Gasket
which induces the trivial element of the group $K^1(C(SG))$. 
This may seem a bit strange, but the idea is that $ZGT$ -- an acronym for {\it Zero Gasket Triple} -- will carry all the geometric information on
the geodesic distance, the Hausdorff dimension, and the volume form, but induce the zero-element in the K-homology group. 
If we then add any
unbounded Fredholm module by direct sum to $ZGT$, 
the $K^1$ element induced by the sum equals that of the added unbounded Fredholm module.

The spectral triple $ZGT$ is defined in
very much the same way as the old gasket triple  
 mentioned in Section \ref{OST}.
It is the direct sum of the unbounded Fredholm 
modules associated to all the upside 
down triangles $\nabla_{m,k}$ plus the unbounded Fredholm module 
coming from the outer triangle $\Delta_{0,1},$ with the
orientation reversed:

\begin{Prop} \label{Sum} 
The direct sum of unbounded Fredholm modules 
\begin{displaymath} 
ZGT:=\overline{UFM(\Delta_{0,1})}\oplus\left(
\oplus_{m=1}^{\infty}(\oplus_{k=1}^{3^{m-1}}{UFM(\nabla_{m,k}}))\right )
\end{displaymath} 
is a spectral triple.  The Hilbert space is
  denoted $H_{ZGT},$ the representation of $C(SG)$ on 
$H_{ZGT} $ and the Dirac operator on $H_{ZGT}$ are named  
$\pi_{ZGT}$ and $D_{ZGT}.$   The bounded Fredholm module coming from the polar
  decomposition of $D_{ZGT}$ induces the trivial element of the
 group $K^1(C(SG))$.  
\end{Prop}  

\begin{proof}
The arguments from \cite{CIL}, pages 27--28, may be copied and  
show that the direct sum is a spectral triple.

The results (ii) and (iii) in Proposition \ref{KUFM}  show that the
corresponding element in $K^1(C(SG))$ is trivial. 
\end{proof}

\begin{Thm}
The spectral triple $ZGT$ has the following
geometric properties:
\begin{itemize}
\item[(i)] The metric induced by the $ZGT$ is the geodesic
  distance. 
\item[(ii)] The $ZGT$ is summable for any positive $s > \log 3/\log 2.$  Its zeta-function $\zeta_{ZGT}(s)$  is meromorphic with a simple pole at  $\log 3/ \log 2$ and it is given by 
$
 \zeta_{ZGT}(s) \, = \, 2 \frac{(2^s-1)(2^s-2)}{2^s -
  3} \cdot \zeta(s).$
\item[(iii)] Let $\mu$ denote the Hausdorff
  probability measure of dimension $\log3 /\log 2$ on the Sierpinski Gasket. Then for any 
Dixmier trace and any continuous function $g$ in $C(SG)$ we have
 \begin{displaymath}
\Tr_{\omega}\left (|D_{ZGT}|^{-\frac{\log 3}{\log 2}}\pi_{ZGT}(g)\right)
 =  \frac{4}{3 \log 3} \cdot \zeta \left ( \frac{\log 3}{\log 2} \right ) \cdot \int_{SG}g(x) d\mu(x). \end{displaymath}
\end{itemize}
\end{Thm}
\begin{proof}
Except for the computation of the zeta-function, the proof here can 
be copied  from \cite{CIL}, Theorem 8.4, Proposition 8.6.
With respect to the zeta-function, we get for any fixed natural
number $l$ that the zeta function for $\Delta_{l,j} $ or
$\nabla_{l,j} $ in a point $s > 1 $  equals 
\begin{align}  \label{zeta1}
\sum_{k \in \bz}2^{-ls}|k + 1/2|^{-s }  &= 2\cdot 2^{-ls} \cdot 2^s\sum_{k \in
  \bn} (2k +1)^{-s}\\ &= 2 \cdot 2^{-ls} \cdot 2^s \cdot (1-2^{-s})\cdot \zeta(s)\\
\label{zeta3} &= 2^{1-ls} \cdot (2^s-1) \cdot \zeta(s)
\end{align}
Since there are $3^{m-1}$ triangles of the form $\nabla_{m,j} $ and
just one of the form $\Delta_{0,1}$ we get that the zeta function for
this spectral triple -- for $ s > \log(3)/ \log(2)$ -- is given by 
\begin{align*}
\zeta_{ZGT}(s) &= 2(2^s-1)\cdot \zeta(s) +
\sum_{m=1}^{\infty}3^{m-1}\cdot 2\cdot 2^{-ms}(2^s -1) \cdot \zeta(s) \\
&= 2(2^s-1)\zeta(s) + 2 \cdot (2^s-1)2^{-s}\frac{1}{1 - (3/2^s)} \cdot \zeta(s)\\
&= 2 \frac{(2^s-1)(2^s-2)}{2^s - 3} \cdot  \zeta (s).
\end{align*}
Further on, by Proposition 4 on page 306 in \cite{Co2}, we obtain 
\begin{align*} 
& \mathrm{Tr}_{\omega}\left (\vert D_{ZGT}\vert ^{- \frac{\log 3}{\log
    2}}\right  )  =  \underset{x \to 1+}{\lim} (x-1){\rm Tr}\left (   \vert D_{ZGT}\vert ^{-x\cdot \frac{\log 3}{\log 2}}\right ) \\ 
  = & \underset{x \to 1+}{\lim} (x-1) \zeta_{ZGT}\left ( x\cdot \frac{\log 3}{\log 2} \right )\\ 
= &  2\underset{x \to 1+}{\lim} (x-1) \frac{\left ( 2^{x \frac{\log 3}{\log 2}}-1\right )\left (2^{x\frac{\log 3}{\log 2}}-2\right )}{2^{x\frac{\log 3}{\log 2}} -
  3} \zeta \left (x \frac{\log 3}{\log 2} \right )\\
 = & \frac{4}{3}\cdot \zeta\left (\frac{\log3}{\log
  2}\right )\cdot \underset{x \to 1+}{\lim} \, \frac{x-1}{3^{x-1}-1}\\
= & \frac{4}{3\log3}\cdot \zeta\left (\frac{\log3}{\log 2}\right ). 
\end{align*}
\end{proof}
\begin{Dfn}\label{UFM}  
For $(f(m,k))_{m,k} \in \Pi_{(m,k) \in HSG}\bz$ let $UFM(f)$ denote
 the direct sum 
\begin{displaymath}
\underset{f(m,k) \neq 0}{\underset{(m,k) \in HSG:}{\oplus}}
\begin{cases} UFM(\nabla_{m,k})\overset{f(m,k)}{\oplus \ldots \oplus} UFM(\nabla_{m,k})  \text{ if } f(m,k) > 0. \\  
\overline{UFM(\nabla_{m,k})}\overset{-f(m,k)}{\oplus \ldots \oplus} \overline{UFM(\nabla_{m,k})}  \text{ if } f(m,k) < 0. \\  
\end{cases}
\end{displaymath}
\end{Dfn}
It is now quite easy to construct a spectral triple which induces any prescribed element in  
the group $K^1(C(SG))$. 
\begin{Dfn} Let $(f(m,k))_{m,k} \in \Pi_{(m,k) \in HSG}\bz$ then $ST(f)$
  denotes the direct sum 
\begin{displaymath} ZGT\oplus
    UFM(f).
\end{displaymath} 
\end{Dfn}
\begin{Thm}  Let $(f(m,k))_{m,k} \in \Pi_{(m,k) \in HSG}\bz.$ Then  
$ST(f)$ is a  spectral triple which induces the geodesic distance on the gasket and the $K^1$-element $(f(m,k))_{m,k}$. 
\end{Thm}

Such a spectral triple will in general not
have the right summability properties.  
In the first place the authors thought that the elements in 
$K^1(C(SG))$ coming from spectral triples that encode the fractal geometry of the gasket ought to be a subgroup of 
$\mathrm{Hom}(K_1(C(SG)), \bz)$ but it is by no means clear or even
true that the direct sum of two arbitrary unbounded  Fredholm modules will make
sense as an unbounded Fredholm module. The problem is analogous to 
the one coming from the addition of unbounded operators, namely that
for two spectral triples $(A_1, H_1, D_1)$ and $(A_2, H_2, D_2)$ associated to a C*-algebra $\mathcal{A}$ the
intersection $A_1 \cap A_2$ may be too small to be a dense 
subalgebra of $\ca.$

This problem does not arise when forming the sum used in the definition of  $ST(f),$ since 
the algebra $A$ is the same for all summands. We infer from the proof of Theorem 8.2 in \cite{CIL} that $A$ is the algebra generated by the real affine functions on the plane, restricted to the gasket. Since $ZGT$ is a summand in any $ST(f)$ it follows that the metric induced by the $ST(f)$ must be the one induced by $ZGT$ which equals the geodesic distance on the gasket.

With respect to summability properties, it follows from the fact that $ZGT$ is a direct summand in $ST(f)$ that $ST(f)$ can only be summable for     
$p > \log 3/\log 2.$ The zeta function for $ST(f),$ say $\zeta_f$,  is easily computed to be  
\begin{align*}
& \zeta_f(s)= \\
& 2 \frac{(2^s-1)(2^s-2)}{2^s - 3} \cdot \zeta(s) +\underset{m=1}{\overset{\infty}{\sum}}  \underset{k=1}{\overset{3^{m-1}}{\sum}} \vert f(m,k) \vert (2 \cdot 2^{-ms} \cdot (2^s-1) \cdot \zeta(s)) \\
& =\zeta_{ZGT}(s)+ 2\cdot (2^s-1) \zeta(s) \left (\underset{m=1}{\overset{\infty}{\sum}}  \underset{k=1}{\overset{3^{m-1}}{\sum}} \vert f(m,k) \vert \cdot 2^{-ms} \right ).
\end{align*}
By applying the root criterion we see that the following result holds.
\begin{Thm}\label{SST}
Let $(f(m,k)) \in \Pi_{HSG}\mathbb{Z}$. Then $ST(f)$ is summable for any $p>\log3 /\log 2$ if and only if
\begin{displaymath}
\underset{m\rightarrow \infty}{\limsup} \left ( \sum_{k=1}^{3^{m-1}}|f(m,k)| \right ) ^{1/m} \leq 3.
\end{displaymath}
\end{Thm}
There is an immediate corollary.
\begin{Cor}
If $ \sup_{m,k}|f(m,k)| <\infty $, then $ST(f)$ is summable for any $p >\log 3/\log 2.$
\end{Cor}
This result implies that the elements of the subgroup of the K-homology group $K^1$ given by bounded sequences are all represented by spectral triples with the right summability properties.  
We have thought of a possible converse to Theorem \ref{SST}, and we conjecture that the following is true.
 \begin{Conjecture}  Let $(A, H, D)$ denote a spectral triple for $C(SG)$
   which induces the geodesic distance on the gasket.  If it is summable for any $p>\log 3/\log 2$ and $(f(m,k)) \in \Pi_{HSG}\mathbb{Z}$ represents its K-homology class then 
\begin{displaymath}
\underset{m \to \infty}{\limsup}\left (\sum_{k=1}^{3^{m-1}}
|f(m,k)|\right )^{1/m} \leq 3.
\end{displaymath}
\end{Conjecture}

At the time of the publication of this article we had no new information on the status of the conjecture. 

 The volume form obtained from a spectral triple $ST(f)$ which is summable for $p >\log 3/\log 2$ clearly majorizes a multiple of the one obtained from the $\log 3/\log 2$-dimensional Hausdorff measure, since $ST(f)$ contains $ZGT$ as a direct summand. On the other hand it is rather obvious that even for a bounded sequence $(f(m,k))$ the values may be unevenly distributed so that the resulting volume form is not proportional to the one coming from the $\log 3/\log 2$-dimensional Hausdorff measure. Below, we discuss conditions on $(f(m,k))$ that assure that the volume form given by $ST(f)$ and the $\log 3/\log 2$-dimensional Hausdorff measure are multiples of each other. We first give a definition, then state our results. 

\begin{Dfn}
The sequence $( f(m,k))$  is said to be boundedly almost invariant, if
 the sequence of  averages 
$$a(m) := \sum_{k=1}^{3^{m-1}} \frac{|f(m,k)|}{3^{m-1}} $$
 is bounded and 
$$ \underset{m \to \infty}{
\limsup}\left ( \sum_{k=1}^{3^{m-1}} \left \vert \, |f(m,k)| - a(m)\right \vert \right )^{1/m} < 3.$$ 
\end{Dfn}

Let $D_{f}$ and $(H_f, \pi_f) $ denote the Dirac operator and the representation of  the spectral triple
$ST(f)$. 

\begin{Lemma}
  Let $(f(m,k))$ be a boundedly almost invariant sequence. Then
 \begin{displaymath}
|D_f|^{-\log 3/\log 2} \in \cl^{(1, \infty )}(H_f),  \quad ( \text{See beginning of \cite{Co2}}, \rm{IV.2.\beta}).
 \end{displaymath} 
\end{Lemma}

\begin{proof}

According to our description of the spectral triple $ST(f)$ and the
statement in front of Theorem 3.7 it follows that the $\z$-function
 for $s > \log3/\log 2 $  is given as the infinite sum of non-negative values

\begin{displaymath}
 \Tr(|D_f|^{-s}) = \zeta_f(s)=\zeta_{ZGT}(s)+ 2(2^s-1) \zeta(s) \left (\underset{m=1}{\overset{\infty}{\sum}}  \underset{k=1}{\overset{3^{m-1}}{\sum}} \vert f(m,k) \vert \cdot 2^{-ms} \right ).
\end{displaymath}

The theory we will use to prove the lemma is all based on limits as
$x \to 1_+$ rather than $s \to \log 3/\log 2,$ so we will replace
the $\z$-function with the positive function $T(x)$ defined for $x>1$ by 
\begin{displaymath}
T(x):= \z_f\left (x \cdot \frac{\log 3}{\log 2} \right).
\end{displaymath} 
Since $(f(m,k))$ is almost boundedly invariant we define the sequence
of averages $(a(m))$ as above and let 
\begin{displaymath}
A := \sup\{a(m): m \in
\bn_0 \} +1.
\end{displaymath}
Then there exists a positive $d<3$ and a natural number
$N_1$ such that
$$ \forall m \geq N_1: \quad \sum_{k=1}^{3^{m-1}}\left |
|f(m,k)|- a(m) \right | < d^m.$$ 
Since all elements in the sum which
defines $\z _f(s)$ are non-negative we may rearrange the sum  as we
please. Further since we are interested in the behaviour of
$(x-1)T(x)$ as $x \to 1_+$ any finite number of terms in the sum
defining $T(x)$ may be left out, when we prove that $(x - 1) T(x) $ is
bounded on the interval $(1, \infty).$ We know already by Theorem 3.3 that
the function $\z_{ZGT}(x\log3/\log 2)$ can be extended to a
meromorphic function with a simple pole at $1.$ Let us then look at
the remaining sum and let us begin to sum for $m \geq N_1.$ To make
the notation easier we will still use the variable  $s = x
\log 3/\log 2$ and we get; 
\begin{align*}
  & 2\cdot (2^s-1) \zeta(s) \sum_{m=N_1}^{\infty}
   \sum_{k=1}^{3^{m-1}} \vert f(m,k) \vert \cdot
    2^{-ms}  \\
  & \leq 2(2^s -1)\z(s) \sum_{m=N_1}^{\infty}2^{-ms}\left[ 3^{m-1}a(m)\,
  + \,  \sum_{k = 1
    }^{3^{m-1}} \left \vert \vert  f(m,k) \vert-a(m) \right \vert  \right ] \\
& \leq 2(2^s -1)\z(s) \sum_{m=N_1}^{\infty}2^{-ms}\left [3^{m-1}A +
 d^m \right ] \\ 
& =  2(2^s  -1)\z(s)\left [\frac{A}{3}\left (\frac{3}{2^s}\right )^{N_1}\frac{2^s}{2^s-3} \, + \,\left (\frac{d}{2^s}\right )^{N_1}\frac{2^s}{2^s-d}\right ].
\end{align*}

From here it follows that the function $\z_f(s)\cdot \left ( s-\frac{\log 3}{\log 2} \right )$
is bounded for $s \in \left ( \frac{\log 3}{\log 2}, \infty \right ),$ and 
then $T(x)(x-1)$ is bounded on the  interval $(1, \infty).$ The lemma
then follows from  \cite{GBVF}, Lemma 7.19 and Lemma 7.20. 
\end{proof}

\begin{Prop}
  Let $(f(m,k))$ be a boundedly almost invariant sequence, with
  corresponding spectral triple $ST(f)$, and $\omega$ an ultrafilter
  on $\bn$ which induces a Dixmier trace on the ideal $\cl^{(1, \infty
    )}(H_f), $ then the functional $\f $ on $C(SG)$ defined for $g $
  in $C(SG)$ by 
  $$\f(g) := \Tr_{\omega}(|D_f|^{-\log 3/\log 2}\pi_f(g))$$ is a
  multiple of the Hausdorff integral on the Sierpinski gasket.
 \end{Prop}

\begin{proof}
  It is well known that the trace property of $\Tr_{\omega} $
  implies that $\f$ becomes a bounded positive linear functional on
  $C(SG).$  We show that the measure assigned by $\f$ to the portion of  Sierpinski gasket which is contained inside or on each $\Delta_{n,j}$ equals $3^{-n}$ of the measure assigned to the entire gasket. This will imply that $\f$ is a multiple of the  Hausdorff integral on the gasket since the $\log 3/ \log 2 $-dimensional Hausdorff measure is the unique measure on the gasket  satisfying the foregoing scaling property (\cite{Ba}).
  First of all we will only need to study the Dixmier trace coming
  from the unbounded Fredholm module $UFM(f),$ so we will introduce
  some notations. We will let $\r_{m,k}$ denote the representation
  consisting of $\vert f(m,k) \vert $ copies of the standard representation of
  $C(\nabla_{m,k}), $ and we will let $K_{m,k}$ denote the
  corresponding Hilbert space. On this space we define $\tilde D_{m,k}
  $ as the corresponding amplified Dirac operator.  Then we define the
  Hilbert space $K,$ the representation $\r$ and the Dirac operator as
  the direct sum of these objects over indices $m,k$ with $f(m,k) \neq
  0.$ By the lemma above $|\tilde D|^{-\log 3/\log 2}$ is in
  $\cl^{(1, \infty)}.$

 Take $F_{n,j}$ to be the triangle $\Delta_{n,j}$ together with its interior. Let us then consider the functions on $SG$ which are the
  characteristic functions $\chi_{1,j}\, \, \, , j \in \{1, 2, 3 \},$ for
  the portion of  Sierpinski gasket which is contained in $F_{1,j}.$ These functions are not
  continuous on the gasket, but we  see that there is a unique
  way to extend the representation $\r$ to be defined on these
  functions too, namely by 
   defining $\r(\chi_{1,j})$ as the orthogonal projection 
 from $K$ onto the subspace $K_j$
  defined by $$K_j := \underset{\nabla_{m,k} \subset
    F_{1,j}}{\oplus} K_{m,k}.$$ 
To compute
  $\rm{Tr}_{\omega}\left (\r(\chi_{1,j})|\tilde D|^{-\log 3/\log 2 } \right ),$
 we have to go through summations
  over decreasing eigenvalues of $|\tilde D|^{-\log 3 /\log 2 }$
 corresponding to the eigenvectors which are 
  contained in the space  $K_j.$ Here each nonzero
  eigenvalue for $|\tilde D|^{-\log 3/\log 2}$ is of the form
  $$3^{-m}|1/2+l|^{-\log 3/\log 2} \text{ for } m \in \bn,
  \quad l \in
  \bz.$$
The multiplicity of such an eigenvalue is $$M(m):=\sum_{\nabla_{m,k} \subseteq
    \Delta_{1,j}} |f(m,k)| .$$ 
In the partial sums used to define the
  Dixmier trace $$\rm{Tr}_{\omega}(\r(\chi_{1,j})|\tilde
  D|^{-\log 3/\log 2}),$$  we will now replace
  the expression $$M(m)3^{-m}|1/2+l|^{-\log 3/\log 2 }$$ by
  $$3^{m-2}a(m)3^{-m}|1/2+l|^{-\log 3/\log 2}, $$ and then we need a
  correction term $$\left ( M(m)-
  3^{m-2}a(m) \right )3^{-m}|1/2+l|^{-\log 3/\log 2}.$$
 
  By checking partial sums of sums corresponding to decreasing
  eigenvalues, one can se that the sums involving $a(m)$ will be
  exactly the same as those used to compute
  $\Tr_{\omega}(|\tilde D|^{-\log 3/\log 2}),$ except that they are
  all scaled by the factor $1/3.$ For a suitable $d<3,$ the sums of the correction terms may
  be dominated as follows:

\begin{align*} & \sum_{l \in \bz} \sum_{m \in \bn} |M(m)- 3^{m-2}a(m)|
  3^{-m}|1/2+l|^{-\log 3/\log 2} \\
 &=  4\zeta \left ( \frac{\log 3}{\log 2}\right )\sum_{m \in \bn} |M(m)- 3^{m-2}a(m)|
  3^{-m}\\
  & \leq  4\zeta \left ( \frac{\log 3}{\log 2}\right )\sum_{m\in \bn} 3^{-m}\sum_{\nabla_{m,k} \subseteq
    \Delta_{1,j}}\left \vert \vert  f(m,k) \vert  - a(m)\right \vert \\
&  \leq 4\zeta \left ( \frac{\log 3}{\log 2}\right )\sum_{m\in \bn}3^{-m} d^m \\ 
& \leq  4\zeta\left ( \frac{\log 3}{\log 2}\right )\frac{d}{3-d}.\end{align*} 

Since the sum of all these correction terms is finite, the Dixmier
trace  $\Tr_{\omega}\left ( \r(\chi_{1,j})|\tilde D|^{-\log 3/\log 2} \right )$
is one third of $\Tr_{\omega}\left (|\tilde D|^{-\log 3/\log 2} \right).$ Repeated use of this argument will show that the measure assigned by $\f$ to the portion of  the Sierpinski gasket which is contained in $F_{n,j}$ equals $3^{-n}$ of the measure assigned to the entire gasket. 
\end{proof} 
\begin{exmp}
  The boundedly almost invariant sequences do not form a
  subgroup of $K^1(SG)$: There exist two boundedly almost
  invariant sequences $f$ and $g$ such that for the sum sequence
  $h(m,k) := f(m,k) + g(m,k) $, the Dixmier trace associated to $ST(h)$
  is not a multiple of the Hausdorff integral. 
\end{exmp}

\begin{proof}
The gasket has the three upward pointing triangles $$\Delta_{1,k},\, k=1,2,3.$$
We let $(f(m,k))$ and $(g(m,k))$ denote the sequences
in $K^1(SG))$ given by $f(m,k) =1$ for all indices
$(m,k)$ and  
\begin{displaymath} 
  g(m,k) = \begin{cases} -1 \text{ if } \nabla_{m,k} \subset
    F_{1,1}\\
    1 \text{ else, }\end{cases} \end{displaymath}
    respectively. Here $F_{1,1}$ is as in the proof of Proposition 3.12. 
Clearly both  are boundedly almost invariant, but we shall
see that their sum $h(m,k) := f(m,k)+ g(m,k)$ is not.
We obtain $a(m)=4/3$ for the sequence of averages  for  $|h(m,k)|$.
Hence 
\begin{align*}
\lefteqn{\underset{m\to \infty}\limsup \left (\underset{k=1}{\overset{3^{m-1}}{\sum}} \left \vert \vert h(m,k) \vert  -\frac{4}{3} \right \vert \right ) ^{1/m}}\\
  &= \underset{m\to \infty}\limsup
\left (\frac{1}{3}\cdot \frac{4}{3}\,3^{m-1} + \frac{2}{3}\cdot \frac{2}{3}\, 3^{m-1}\right )^{1/m} \\ & =
\underset{m\to \infty}\limsup \left ( 8\cdot 3^{m-3} \right ) ^{1/m} = 3,
\end{align*} 
so the sequence $(h(m,k)) $ is not boundedly almost invariant.\\
\indent For any continuous positive function with support inside $F_{1,1}$, the Dixmier trace  induced by $ST(h)$ equals 
that of the Dixmier trace coming from ${ZGT}$. This will
not be the case, however,  for a positive continuous function supported inside 
$F_{1,2}.$ Hence the Dixmier trace corresponding to $ST(h)$ is not
a multiple of the Hausdorff integral.
\end{proof}
The example shows that in order to obtain a subgroup of sequences
$(f(m,k))$ such that the Dixmier trace generated by $ST(f)$ is
 a multiple of the Hausdorff integral, we can not do with bounded sequences. On the
other hand there is a subgroup where all elements behave nicely.

\begin{Dfn} 
The group $c_1(HSG,\bz)$ is defined as the set of sequences  $(f(m,k))$ such
that
\begin{displaymath}
\exists t \in \bz\, \exists M \in \bn\, \forall m \geq M \, \forall k
\in \{1, \dots, 3^{m-1}\}:\,f(m,k) = t.\end{displaymath}
\end{Dfn}
\begin{Prop} 
For $(f(m,k)) $  in $c_1(HSG, \bz),$ the 
Dixmier trace associated to the spectral triple $ST(f)$ is a multiple
of the Hausdorff integral.
\end{Prop}

\begin{proof}
The sequences in $c_1(HSG, \bz)$ are all boundedly almost invariant.
\end{proof}

\section{Aspects of Minimality of $ZGT$}

In \cite{Co3} Connes formulates seven axioms for a spectral triple
\newline  $(\mathcal{A}, \mathcal{H}, D)$ with $\mathcal{A}$ a commutative algebra.
In \cite{Co5} he shows that five of these axioms (in a slightly stronger form) suffice in order to characterize the spectral triples associated to smooth compact manifolds, i.e. from these axioms one can construct a smooth oriented compact (spin$^c$) manifold $X.$ 

The fractal we study is by nature far away from being a smooth manifold, and this is why we can not expect to construct  spectral triples for the Sierpinski gasket which satisfy all of Connes' axioms.
%The fractals we study here are by nature far away from being smooth
%manifolds and we think that this 
%is the reason why we have not been able to
%apply Connes' results to our present investigation.
One of the problems our spectral triples raise is that
our algebra of {\em differentiable functions (Lipschitz functions), } according to
\cite{CIL}, Theorem 8.2, is the algebra of functions on the gasket generated by
the restrictions of the affine functions. As we remark in that paper, an affine
function restricted to a triangle gives a continuous function with
constant slopes along the edges. It then belongs to the domain of the
{\em Dirac} operator, but its derivative will not, unless the
function is constant. It seems
hopeless to get a regularity condition, i.e. to talk about smooth functions in this case.
%as the one Connes asks for \cite{Co3, Co5}, to be in fulfilled in our case. 
We will therefore stay inside the
set-up where we consider spectral triples defined as direct sums of
of unbounded Fredholm modules associated to triangles of the form
$\nabla_{m,j} $ or $\Delta_{0,1}.$ 
Our main result in this section shows that in this set-up, the geometry, i.e. the representation (\cite{Co3}), of the spectral triple $ZGT$
is the minimal one among those which induce the geodesic
distance on the Sierpinki gasket.  Recall that the $ZGT$ spectral triple is the direct sum of
the unbounded Fredholm modules associated to the big outer triangle
and all those associated to the upside down triangles. As the gasket is the closure of the
increasing sequence of graphs $G_n$ in Figure \ref{SG2}, it is
clear that one can leave out the big outer triangle and any
finite number of upside down triangles and still gets the gasket as the
closure. The corresponding restricted sum of unbounded Fredholm
modules will then induce a faithful representation and hence a
spectral triple, but the geodesic distance will not be recovered.\\
\indent We consider the case of a spectral triple which is obtained from $ZGT$ 
by leaving out one of the triangles. Let us denote such a spectral triple by $RGT$
({\em reduced gasket triple}). We will show:
\begin{Prop} \label{d0}
The metric $d_{RGT}$ induced by  any $RGT$ is equivalent to the geodesic distance $d_g$ on the gasket  but it does not coincide with it: There exist points $x,y,$ such that  
\begin{equation}
d_{RGT}(x,y)\geq \frac{3}{2}d_g(x,y).
\end{equation} 
\end{Prop}
The rest of this section is devoted to a proof of this proposition 
divided into several  small arguments.
A slight modification of the proof of Theorem 8.13 in \cite{CIL}
shows that for any spectral triple of the form $ST(f)$ the metric
induced by the distance formula applied to this triple is $d_g.$
It is well known, see for instance \cite{Ba}, that for points $x, y$ on
the gasket 
\begin{displaymath}
\Vert x - y \Vert \leq d_g(x,y) \leq 8 \Vert x-y \Vert ,
\end{displaymath}
so that the geodesic metric is equivalent to the Euclidean metric.\\
\indent To enter into the proof of the proposition, let us look at the triple
$ZGT$ and suppose that we take out the summand corresponding to the
outer triangle $\Delta_{0,1}.$ It will become clear that our
considerations may be transferred to the situation where we
instead remove one of the triangles $\nabla_{m,k}.$ In order to
clarify things we will give some definitions and figures. Define the graphs
\begin{displaymath} H_n \, :=\, \underset{m\leq n}{\cup} \underset{k\leq 3^{m-1}}{\cup} \nabla_{m,k},\, n=1,2, \ldots, \infty.
\end{displaymath}
See the illustration below for a drawing presenting  $H_3$ and a part of $H_4.$ 

\newcommand{\drawtriangles}{
\draw (0,0) -- (1,1.732) -- (-1,1.732) -- cycle;
\draw[scale=.5,xshift=-2cm] (0,0) -- (1,1.732) -- (-1,1.732) -- cycle;
\draw[scale=.5,xshift=2cm] (0,0) -- (1,1.732) -- (-1,1.732) -- cycle;
\draw[scale=.5,yshift=3.464cm] (0,0) -- (1,1.732) -- (-1,1.732) -- cycle;
}
\begin{figure}[H]
\begin{tikzpicture}
\drawtriangles
\begin{scope}[xshift=4cm] \drawtriangles \end{scope}
\begin{scope}[xshift=2cm,yshift=3.464cm] \drawtriangles \end{scope}
\draw (2,0) -- (4,3.464)  --
(0,3.464) node[anchor=east] {$Q$} -- cycle; 
\draw (.5,3.464) -- (.75,3.897) -- (.25,3.897) -- cycle;
\draw[xshift=.5cm,yshift=.866cm] (.5,3.464) -- (.75,3.897) -- (.25,3.897) --cycle;
%\draw (1,3.464) node[anchor=north] {$T$};
%\draw (1.5,4.33) node[anchor=west] {$S$};
%\draw (.5,4.33) node[anchor=east] {$R$};
\draw (0.7,5) node[anchor=south] {$P_2$};
\draw (2,-1) node [anchor=west] {$H_3$ and part of $H_4$};
%\draw (1.9,4.33) node[anchor=west] {$W$};
\draw (-1.7,1.8) node[anchor=west] {$P_1$};
%\draw[dashed](-1,1.732)--(0,3.464);
\draw (1,1.8) node[anchor=west] {$R$};
\draw (1.7, 3.2) node[anchor=west] {$S$};
%\draw (4, 3.464) node[anchor=west] {$T$};
%\draw (2,0) node[anchor=north] {$V$};
\draw[dashed](-2,0)--(2,6.93)--(6,0)--cycle;
\draw[xshift=-.5cm,yshift=-.866cm] (.5,3.464) -- (.75,3.897) -- (.25,3.897)
--cycle; 
\draw[xshift=-1cm,yshift=-1.732cm] (.5,3.464) -- (.75,3.897) -- (.25,3.897)
--cycle; 

\end{tikzpicture}
\caption{\label{ZG1}}
\end{figure}
On each graph $H_n$ a metric $d_n$ is given by the geodesic distance on that graph.  As $H_n \subseteq H_{n+1}$ the sequence $d_{n+k}(x,y) $ is
decreasing for any points $x,y$ from $H_n$ as $k \to \infty,$ and we obtain a metric on
$H_{\infty} $ by  
\begin{displaymath}
d_{\infty}(x,y):= \underset{k \to \infty}{\lim} d_{n+k}(x,y).
\end{displaymath}
\begin{Lemma}
For any points $x, y$ in $H_{\infty}$
\begin{displaymath}
d_g(x,y) \, \leq \, d_{\infty}(x,y) \leq 2 d_g(x,y).
\end{displaymath}
\end{Lemma}

 \begin{proof} 
Fix $n$ such that $x,y \in H_n.$ It is known from the analysis in \cite{CIL}, see  Lemma 8.11 and the following remarks, that a shortest path from $x$ to $y$
can be found along the edges of $G_n=H_n \cup \Delta_{0,1}.$   As any piece of this path going through $\Delta_{0,1}$ can be replaced by two pieces of the same length in $H_n,$ we have 
\begin{displaymath}
d_g(x,y) \, \leq \, d_{\infty}(x,y) \leq d_n(x,y) \leq 2 d_g(x,y). \qedhere
\end{displaymath}
\end{proof}

Since the set of points in 
$H_{\infty}$ is dense in the gasket, it is possible to extend $d_{\infty}$ by continuity to a metric $d_0$
on the gasket, and we have:
\begin{Cor} For any pair of points $x, y$ on the Sierpinski gasket
  \begin{displaymath}
d_g(x,y) \leq d_0(x,y) \leq 2 d_g(x,y),
\end{displaymath} 
so that the metric $d_0$ is equivalent to the geodesic distance.
\end{Cor} 
We will show next that $d_0$ does not  coincide with the geodesic distance. As mentioned above  
let $RGT$ denote the {\em reduced } spectral triple obtained as
the direct sum of all the unbounded Fredholm modules associated to all
the triangles $\nabla_{m,j} $ (thus leaving out $\Delta_{0,1}$),  and let $D_{RGT}, \, \pi_{RGT}$ be the
corresponding Dirac operator and representation.  
We consider two vertices $P_1$ and $P_2$ on $\Delta_{0,1}$ as sketched in Figure \ref{ZG1} and shall show that
\begin{displaymath}
d_n(P_1,P_2)=\frac{3}{2}d_g(P_1,P_2), \, \text{for } n\geq 3.
\end{displaymath}
This implies that 
\begin{equation}
d_0(P_1,P_2)=\frac{3}{2}d_g(P_1,P_2).
\end{equation}
We first observe that the function $ d: SG \to
\br_+\cup \{0\}$ defined by $d(x ) := d_0(P_1,x)$ satisfies 
\begin{displaymath}
\Vert [D_{RGT}, \pi_{RGT}(d)] \Vert  = 1,
\end{displaymath}
 and this follows as in the proof of \cite{CIL},
Lemma 8.12, noting that on
any edge in $H_n$ the function $g$ is differentiable except at at most
one point and the numerical value of the derivative, when defined,
is bounded by 1. Hence
\begin{equation}
d_{RGT}(P_1,x) \geq d_0(P_1,x).
\end{equation}
We consider $d$ on the edge of $\Delta_{0,1}$ which contains $P_1$ and $P_2$ and abbreviate $\alpha =d_g(P_1,P_2).$ 
Clearly any shortest path from $P_1$ to $P_2$ will pass through $Q$ and 
\begin{displaymath}
d_2(P_1,P_2)=2\alpha , \, d_3(P_1,P_2)=\frac{3}{2}\alpha. 
\end{displaymath}
Next, considering the triangles with vertices $P_1, R,Q$ and $Q,S,P_2$ we see that determining a shortest path in $H_4$ is analogous to determining a shortest path
in $H_3$ by scaling,  and $d_4(P_1,P_2)=\frac{3}{2}\alpha.$ Iteration then shows that 
\begin{displaymath}
d_n(P_1,P_2)=\frac{3}{2}\alpha, \, \text{for all } n\geq 3.
\end{displaymath}  
Taking the limit as $n\rightarrow \infty$ we obtain $(2)$  and, in view of $(3),$ inequality $(1)$ in Proposition 4.1
\begin{Co}
 Proposition \ref{d0}  above shows that no summand in our
  spectral triple $ZGT$ can be left out, if the geodesic distance shall
  be obtained via Connes' formula. On the other hand one may leave out any
  finite number of summands and still get a spectral triple, because
  the union of the remaining triangles $\nabla_{m,j} $ is dense
  in the gasket. Although such a spectral triple will not reproduce the
  geodesic distance, the volume form will still be
 the one obtained from  $ZGT.$ 
 \end{Co}
\section{The Sierpinski Pyramid:Constructions, K-theory and K-homology}
The Sierpinski pyramid is the 3-dimensional version of the Sierpinski gasket.   One starts, for example, with a solid regular tetrahedron $P_0$ (figure below on the left) and divides it into eight identical regular tetrahedra. One then cuts out all the smaller tetrahedra except for the ones at the vertices of the starting tetrahedron, obtaining $P_1$ (figure below on the right). 
\begin{figure}[H]  
\includegraphics[scale=0.14]{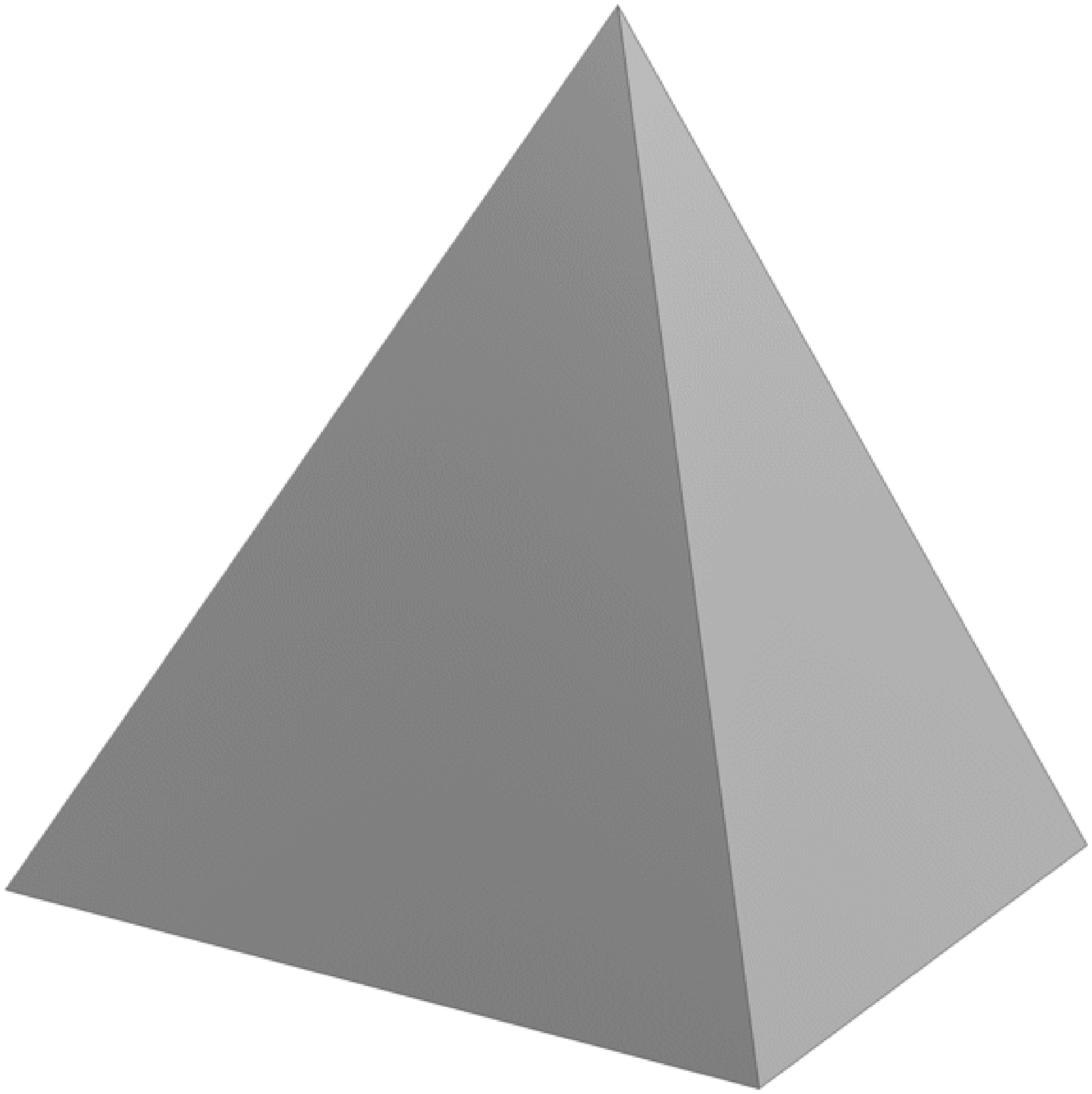}   
\includegraphics[scale=0.14]{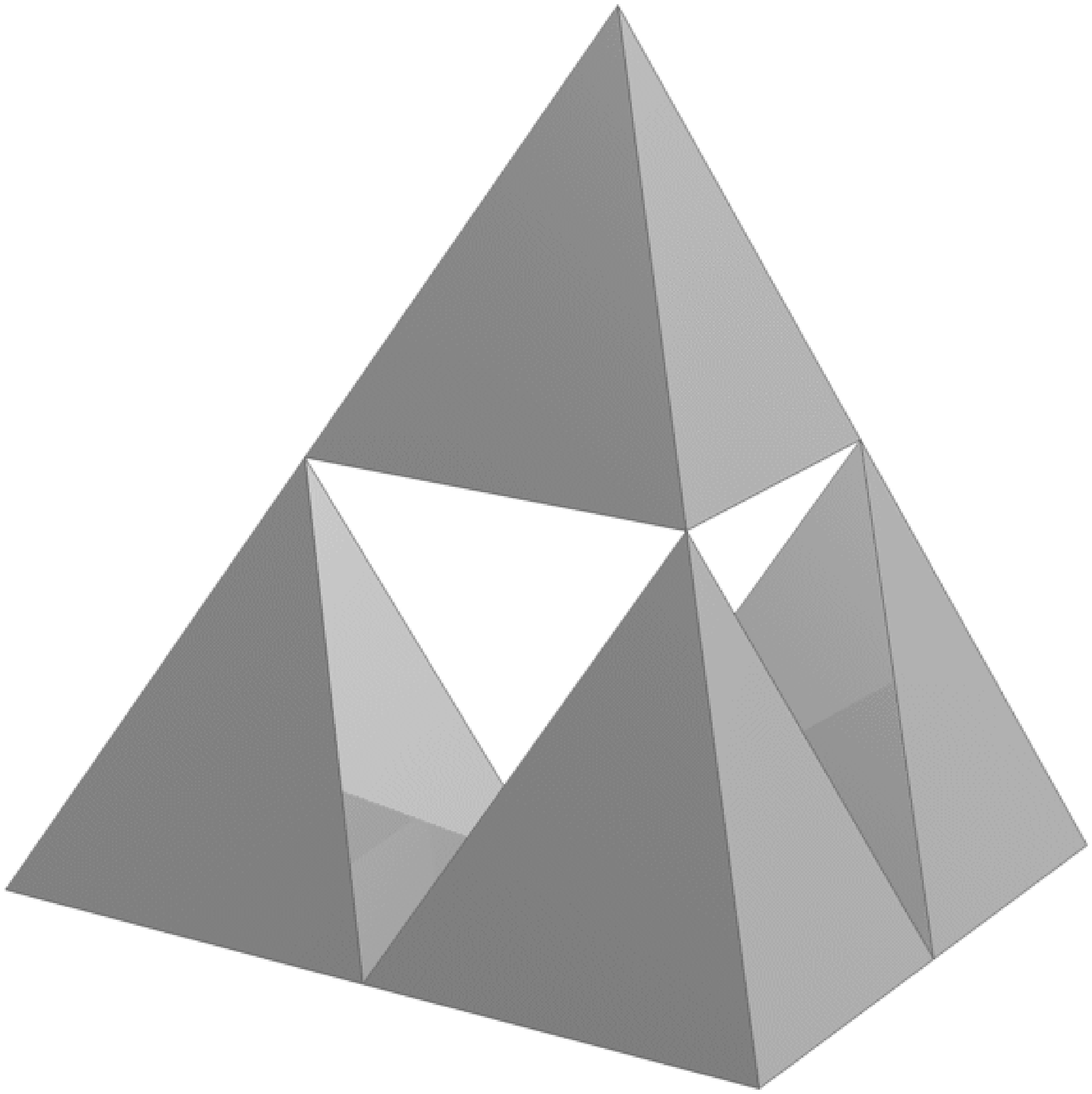}
\caption{\label{SP1}}   
\end{figure}
One iterates this procedure (see below the illustration of $P_2$).
 \begin{figure}[H]   \includegraphics[scale=0.14]{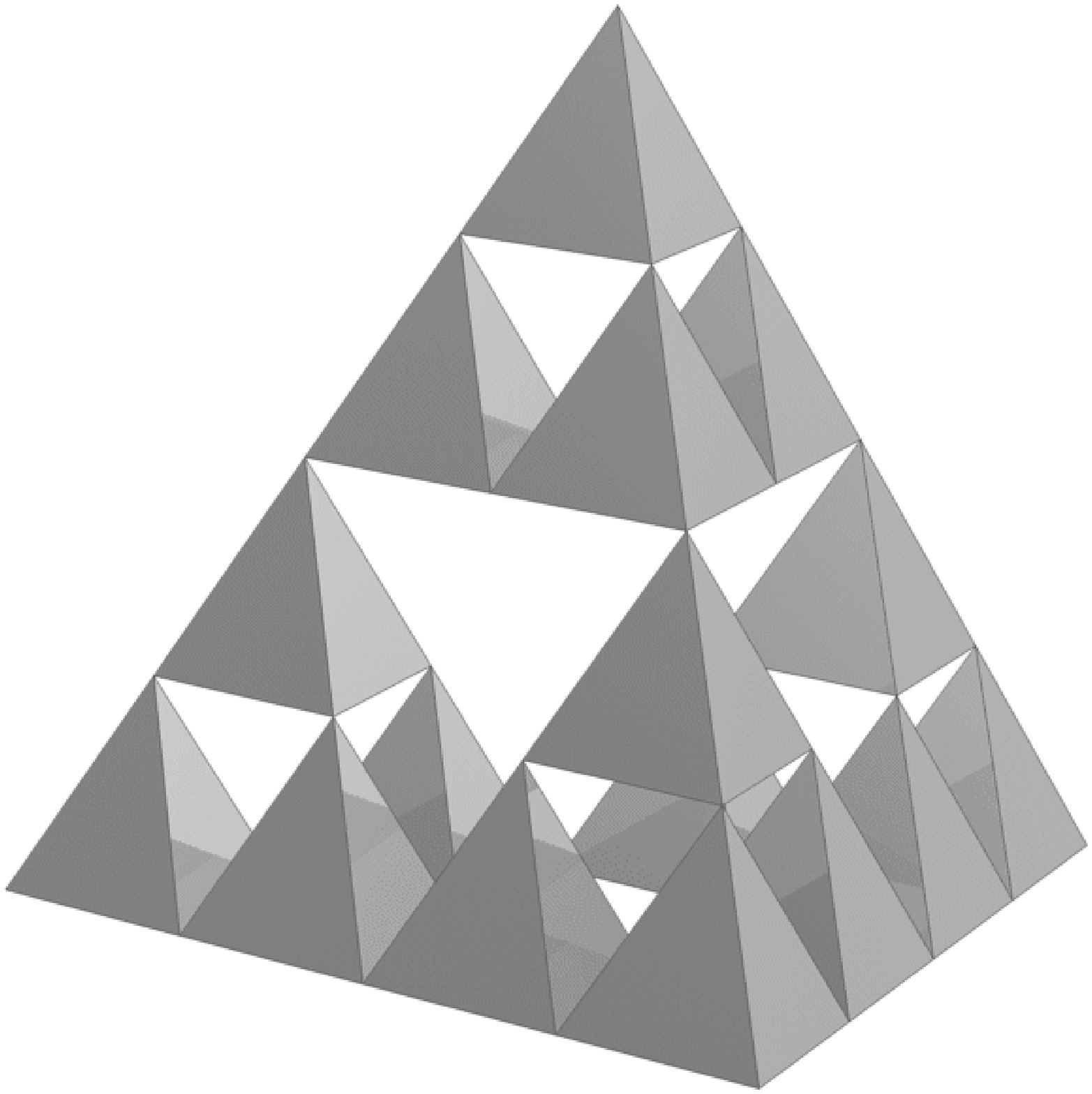} 
 \caption{\label{SP2}}   
\end{figure}
(The drawings are  done by J\"urgen  Meier, www.3d-meier.de.)  \\
\indent In the $n$-th step, one cuts away $4\cdot 4^{n-1}$ regular tetrahedra, or $4^{n-1}$ regular octahedra,  with side length equal to $2^{-n}$ of that of the original one. The Sierpinski pyramid $SP$ is then the limit of this decreasing sequence of compact subsets $P_n$, of $P_0$.
For any $n\in \mathbb{N}, P_n$ is the union of $4^n$ solid tetrahedra, say $P_{n,k}, 1\leq k \leq 4^n$ with the side length  $2^{-n}$ of that of $P_0$. One calculates easily the Hausdorff dimension of this fractal.  Indeed the Sierpinski pyramid is a self-similar set satisfying the open set condition and can be constructed out of 4 similarities of ratio $1/2$. This implies that its Hausdorff dimension $s$ is the solution of the equation $4\cdot(1/2)^s=1$ and thus  $s=2$.
The C*-algebra of the continuous functions on the Sierpinski pyramid $C(SP)$ is the direct limit of C*-algebras $C(P_n)$ and the homomorphisms given by the restriction maps. We will use this picture to compute $K_1(C(SP))$.  We start by considering $P_1$. If we retract each $P_{1,k}, 1\leq k \leq 4$ along oblique edges, as shown below,
\begin{figure}[H]
 \includegraphics[scale=0.15]{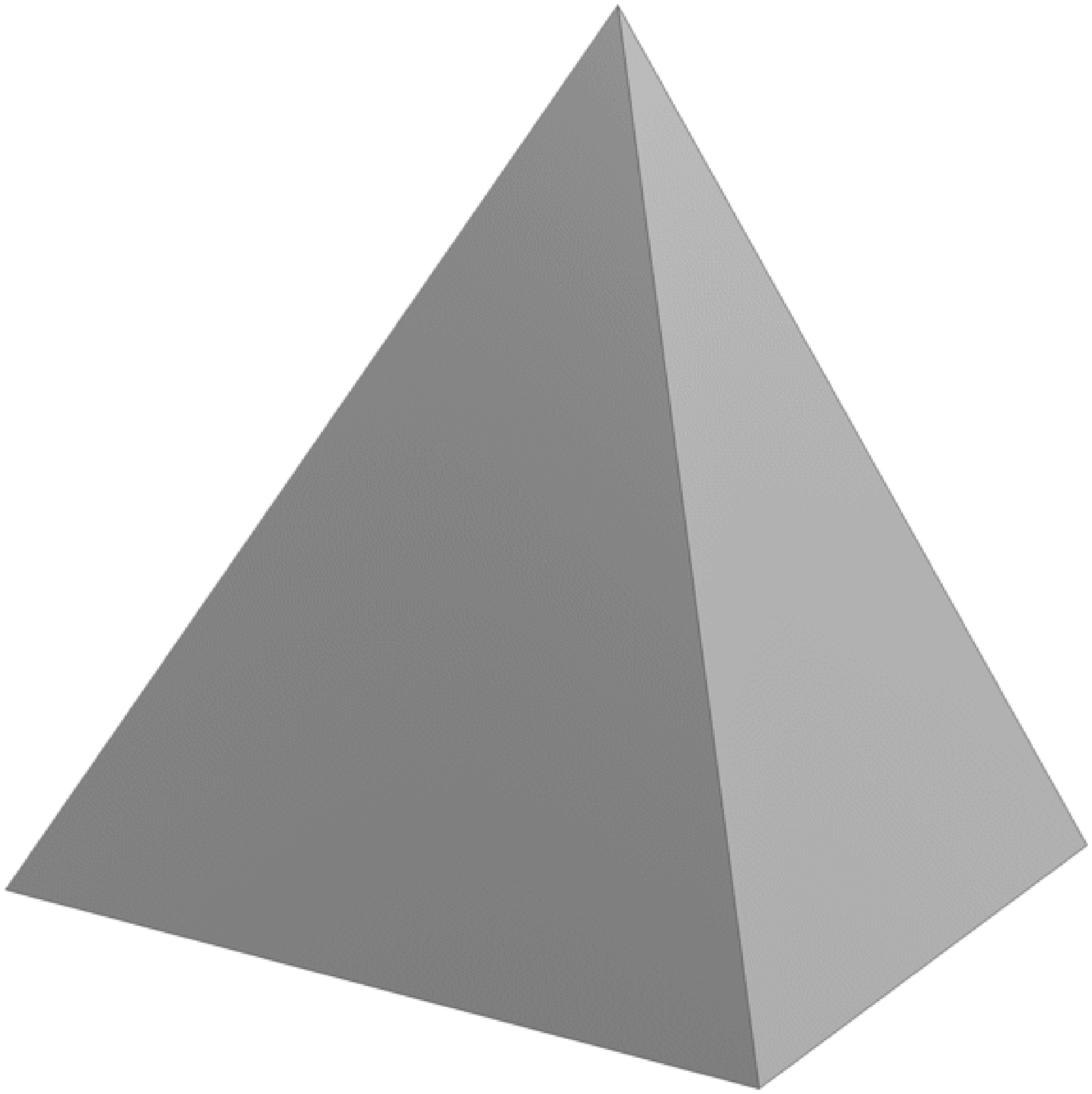} 
  \begin{tikzpicture}[scale=.36]
\draw (5.1,5.1) -- (6.2,-3.0);
\draw (5.1,5.1) -- (8.7, -1.2);
\draw (5.1,5.1) -- (0.5,-1.5);
\foreach \x in {(6.2,-3.0), (8.7,-1.2), (0.5,-1.5), (5.1,5.1)}
  \filldraw \x circle (2pt);
\end{tikzpicture}
\caption{\label{SP3}}
   \end{figure}
\noindent we obtain from $P_1$ the graph, below, which is homotopic to the three-leaved rose $R_3$ (recall that every finite (connected) graph is homotopic to a $n$-leaved rose, where $n$ is the number of edges not belonging to a maximal tree of the graph, see, for example, \cite{AB}).
\begin{figure}[H]
\begin{tikzpicture}[scale=1.2]
\useasboundingbox (-2.5,-2.8) rectangle (2,.3);
  \newcommand{\drawfront}[2]
{
    \draw[white,line width=3pt,opacity=1.0]  (#1) -- (#2);
    \draw (#1) -- (#2);
}
\newcommand{\drawcrown}{
\draw (0,0) -- (-1,-1) -- (-2,-2);
\draw (0,0) -- (0,-1) -- (0,-2);
\draw (0,0) -- (1,-1) -- (2,-2);
\draw (0,0) -- (-1,-1) -- (-2,-2);
\draw (0,0) -- (0,-1) -- (0,-2);
\draw (0,0) -- (1,-1) -- (2,-2);
\draw (-1,-1) -- (-0.8,-2.2);
\draw (0,-1) -- (-0.8,-2.2);
\draw (1,-1) -- (0.8,-2.2);
\draw (0,-1) -- (0.8,-2.2);
\drawfront{-1,-1}{-0.5,-1.5};
\drawfront{1,-1}{-0.5,-1.5};
}
\newcommand{\drawdots}{
\filldraw (-2,-2) circle (1pt);
\filldraw (0,-2) circle (1pt);
\filldraw (1,-1) circle (1pt);
\filldraw (2,-2) circle (1pt);
\filldraw (-0.5,-1.5) circle (1pt);
\filldraw (-0.8,-2.2) circle (1pt);
\filldraw (0.8,-2.2) circle (1pt);
\filldraw (0.8,-2.2) circle (1pt);
\filldraw (0,-1) circle (1pt);
\filldraw (-1,-1) circle (1pt);
\filldraw (0,0) circle (1pt);
}
\newcommand{\drawsmalldots}{
\filldraw (0,0) circle (.5pt);
\filldraw (-1,-1) circle (.5pt);
\filldraw (0,-1) circle (.5pt);
\filldraw (-2,-2) circle (.5pt);
\filldraw (0,-2) circle (.5pt);
\filldraw (1,-1) circle (.5pt);
\filldraw (2,-2) circle (.5pt);
\filldraw (-0.5,-1.5) circle (.5pt);
\filldraw (-0.8,-2.2) circle (.5pt);
\filldraw (0.8,-2.2) circle (.5pt);
}
\draw %Rechteck mit Beschriftungen
(0.1,0.2) node[anchor=east] {$P$};
\draw %Rechteck mit Beschriftungen
(-0.1,-1) node[anchor=east] {$Q$};
\draw %Rechteck mit Beschriftungen
(-1.1,-1) node[anchor=east] {$M$};
\draw %Rechteck mit Beschriftungen
(1.6,-1) node[anchor=east] {$N$};
\draw %Rechteck mit Beschriftungen
(-2.1,-2) node[anchor=east] {$A$};
\draw %Rechteck mit Beschriftungen
(-0.1,-2) node[anchor=east] {$B$};
\draw %Rechteck mit Beschriftungen
(2.5,-2) node[anchor=east] {$C$};
\drawcrown;
\drawdots;
\end{tikzpicture}
% \newpage
\begin{tikzpicture}[scale=3.5]
\useasboundingbox (-1,-.5) rectangle (.6,.4);
\draw[rotate=30] (0,0) .. controls (1,0) and (0,1) .. (0,0);
\draw[rotate=30] (0,0) .. controls (1,0) and (0,-1) .. (0,0);
\draw[rotate=30] (0,0) .. controls (-1,0) and (0,1) .. (0,0);
\filldraw (0,0) circle (0.5pt);
\draw %Rechteck mit Beschriftungen
(0.05,-0.1) node[anchor=east] {$Q$};
%\draw[rotate=30] (0,0) .. controls (-1,0) and (0,-1) .. (0,0);
\draw (0.05,-0.3) node[anchor=east] {$R_3$};
\end{tikzpicture}
\vspace{-.5cm}
\caption{\label{SP5}}
\end{figure}
The three leaves correspond to the holes in the oblique faces of the initial tetrahedron. In the subsequent steps, for any $n \geq 2$  the cutting of an octahedron from one of the tetrahedra of $P_{n-1}$ will produce each time three holes in the oblique faces of $P_{n-1}.$  Hence $P_n$ will be homotopic to a rose having as many leaves as one counts holes in every oblique face of any of the pyramids at the $n$-step in the construction of $SP$, i.e. 
\begin{displaymath}
3+4\cdot 3+\ldots +4^{n-1}\cdot 3 =4^n-1
\end{displaymath}
leaves. The K-theory of the $n$-leaved rose $R_n$ is easily computed using a six term exact sequence and it follows from \cite{RLL}, page 232 that  
\begin{displaymath}
K_0(C(R_n))=\mathbb{Z} \quad \text{and} \quad K_1(C(R_n))=\mathbb{Z}^n.
\end{displaymath}
The singular homology of the $n$-leaved rose is given by $H^1(R_n)=\mathbb{Z}^n$ (\cite{AB,Mu}). We thus conclude 
\begin{displaymath}
K_1(C(P_n))=H^1(P_n)=\mathbb{Z}^{4^n-1}.
\end{displaymath} 
As the direct limit of $K_1(C(P_n))$, $K_1(C(SP))$ will count the holes in every oblique face of any of the small pyramids that arise in the construction of the Sierpinski pyramid\begin{footnote} {This result was obtained in the diploma thesis \cite{Ha} of Stefan Hasselmann supervised by C. Ivan} \end{footnote}, hence 
\begin{displaymath}
K_1(C(SP))=\oplus_{n=1}^{\infty}\mathbb{Z}^{3 \cdot 4^{n-1}}. 
\end{displaymath}
We shall remark here that also for the Sierpinski pyramid 
\begin{displaymath}
K_1(C(SP)) \cong \check{H}^1(SP);
\end{displaymath}
in fact, since $SP$ is the intersection of the $P_n,$  we obtain from \cite{Mu}, Theorem 73.4, that  
\begin{align*}
& \check{H}^1(SP)=\underrightarrow{\lim} H^1(P_n) \quad \text{and} \\
& H^1(P_n)=K_1(C(P_n)). 
\end{align*}
Let $X$ be a compact subset of the space $\mathbb{R}^3$  and $u$ a unitary in some matrix algebra $M_k(C(X)),$ i.e. $u\in U_k(C(X)).$ Then the $K_1$ class $[u]$ of $u$ is not in general represented by a unitary $v$ in $C(X),$ but the above analysis shows that for the Sierpinski pyramid, this is the case. In an unpublished article \cite{Br} (see \cite{BDF}, \cite{Do})  Larry Brown has extended the result in \cite{BDF} to $\mathbb{R}^3,$ i.e. he  showed that  the index map (in Kasparov's picture)
\begin{displaymath}
K^1(C(X)) \ni [(H (\pi),F)] \mapsto \Phi_{(H,F)} \in \text{Hom}(K_1(C(X)), \mathbb{Z})
\end{displaymath}
defined by : 
\begin{displaymath}
\Phi_{(H,F)}([u]) =-\text{Index }(P_k \pi _k(u) P_k), \, 
u\in U_k(C(X)) , \, k \in \mathbb{N} 
\end{displaymath}
is an isomorphism for any compact $X \subset \mathbb{R}^3.$ \\
\indent For the Sierpinski pyramid we conclude 
\begin{displaymath}
\Phi_{(H,F)}([u]) =-\text{Index }(P \pi(u) P), \, 
 u\in U(C(X)). 
\end{displaymath}

\section{A family of Spectral Triples that Generates $K^1(C(SP))$}

We suppose that the initial solid tetrahedron $P_0$ is placed such that one vertex is pointing upwards (as in Figure \ref{SP1}) and has the side length $2\pi/3$. The four equilateral triangles of perimeter $2\pi$ which are the boundaries of the four faces of $P_0$ are denoted $\Delta_{0,k}, \, k\in \{ 1,\ldots,4\}$. 
\begin{figure}[H]
\begin{tikzpicture}[scale=.4]
\draw (0,0) -- (7.4,-2.8) -- (7.4,2.8) -- (5.1,5.1) -- (0,0);
\draw[dashed] (0,0) -- (7.4,2.8);
\draw (5.1,5.1) -- (7.4,-2.8);
\foreach \x in {(0,0), (7.4,-2.8), (7.4,2.8), (5.1,5.1)}
  \filldraw \x circle (1pt);
  
\draw[->, yshift=.3cm] (2.96,-1.12) -- (4.44,-1.68);
\draw[<-, yshift=-.3cm] (2.96,-1.12) -- (4.44,-1.68);

\draw[->, xshift=.3cm] (6.02,1.94) -- (6.48,0.36);
\draw[<-, xshift=-.3cm] (6.02,1.94) -- (6.48,0.36);

\draw[->, xshift=-.2cm, yshift=.2cm] (2.04,2.04) -- (3.06,3.06);
\draw[<-, xshift=.2cm, yshift=-.2cm] (2.04,2.04) -- (3.06,3.06);

\draw[->, xshift=-.25cm] (7.4,0.3) -- (7.4,1.5);
\draw[<-, xshift=.25cm] (7.4,0.3) -- (7.4,1.5);

\draw[->, xshift=.3cm] (5.9,4.3) -- (6.48,3.72);
\draw[<-, yshift=-.3cm] (5.9,4.3) -- (6.48,3.72);

\draw[->, yshift=.3cm] (2.96,1.12) -- (4.44,1.68);
\draw[<-, yshift=-.3cm] (2.96,1.12) -- (4.44,1.68);

\end{tikzpicture}
\caption{\label{SP7}}
\end{figure}

\begin{numbering}\label{numb} 
With $\Delta_{1,1,k}, \, k \in \{ 1,\ldots,8\},$ we denote the eight equilateral triangles of perimeter $\pi$ which are the boundaries of the faces of the cut out octahedron in the first step of the construction procedure of the Sierpinski pyramid. We will consider them numbered as follows 
\begin{enumerate}
\item $\Delta_{1,1,k}, \, k \in \{ 1,2,3\}$ for the boundaries of the holes into an oblique face of $P_0$ generated by cutting out the octahedron, respectively 
\item $\Delta_{1,1,4}$ for the boundary of the hole into the horizontal face of $P_0$, and 
\item $\Delta_{1,1,k}, \, k \in \{ 5,\ldots,8\}$ for the boundaries of the remaining faces of the cut out octahedron.
\end{enumerate}
\begin{figure}[H]
\begin{tikzpicture}[scale=.25]
\draw [dashed] (7.4,-2.8) -- (7.4,2.8) -- (5.1,5.1);
\draw (5.1,5.1) -- (7.4,-2.8);
\foreach \x in {(7.4,-2.8), (7.4,2.8), (5.1,5.1)}
  \filldraw \x circle (2pt);
  
\draw[<-,xshift=-0.24cm,line width=.3pt] (7.4,-0.49)--(7.4,0.49); %AF
\draw[->,xshift=0.24cm,line width=.3pt] (7.4,-0.49)--(7.4,0.49);
  
\draw[->,xshift=0.24cm,line width=0.3pt] (5.77,2.82)--(6.04,1.87); %AC
\draw[<-,xshift=0.24cm,line width=0.3pt] (6.24,1.23)--(6.51,0.28);

\draw[->,yshift=0.17cm,xshift=0.17cm,line width=0.3pt]
(6.57,3.63)--(5.93,4.27);%CF
\draw[<-,yshift=-0.17cm,xshift=-0.17cm,line width=0.3pt]
(6.57,3.63)--(5.93,4.27);

\begin{scope}[xshift=7.4cm,yshift=-2.8cm]
%\draw (7.4,2.8)--(7.4,-2.8) -- (5.1,5.1);
\draw (0,0) --(7.4,2.8)--(5.1,5.1);
\draw (5.1,5.1) -- (0,0);
\foreach \x in {(0,0), (7.4,2.8), (5.1,5.1)}
 \filldraw \x circle (2pt);
 
\draw[<-,yshift=0.24cm,line width=0.3pt] (2.44,0.92)--(3.48,1.32); %AD
\draw[->,yshift=0.24cm,line width=0.3pt] (3.92,1.48)--(4.96,1.88);

\draw[<-,yshift=0.17cm,xshift=0.17cm,line width=0.3pt]
(6.57,3.63)--(5.93,4.27);%BD
\draw[->,yshift=-0.17cm,xshift=-0.17cm,line width=0.3pt]
(6.57,3.63)--(5.93,4.27);

\draw[->,yshift=0.17cm,xshift=-0.17cm,line width=0.3pt] (2.1,2.1)--(3,3);%AB 
\draw[<-,yshift=-0.17cm,xshift=0.17cm,line width=0.3pt] (2.1,2.1)--(3,3);

\end{scope}

\begin{scope}[xshift=7.4cm,yshift=2.8cm]
%\draw (0,0) -- (4.05,-1.55);
\draw[dashed] (0,0)--(4.05,-1.55) -- (7.4,-2.8);
%\draw (7.4,-2.8) -- (7.4,2.8) -- (5.1,5.1);
\draw[dashed] (5.1,5.1) -- (0,0);
\draw (5.1,5.1) -- (7.4,-2.8);

\foreach \x in {(0,0), (7.4,-2.8), (5.1,5.1)}
  \filldraw \x circle (2pt);

\draw[->,yshift=-0.3cm,line width=0.3pt] (1.33,-0.5)--(2.37,-0.9); %DF
\draw[<-,yshift=0.3cm,line width=0.3pt] (2.81,-1.06)--(3.85,-1.47);

\draw[<-,yshift=0.17cm,xshift=-0.17cm,line width=0.3pt] (2.1,2.1)--(3,3);%EF
\draw[->,yshift=-0.17cm,xshift=0.17cm,line width=0.3pt] (2.1,2.1)--(3,3);

\draw[->,xshift=-0.24cm,line width=0.3pt] (6.04,1.87)--(6.31,0.92); %DE
\draw[<-,xshift=-0.24cm,line width=0.3pt] (6.51,0.28)--(6.78,-0.67);

\end{scope}

\begin{scope}[xshift=5.1cm,yshift=5.1cm]
%\draw (0,0) -- (7.4,-2.8) -- (7.4,2.8) -- (0,0);
\draw (0,0)-- (7.4,-2.8);
\draw (0,0)-- (7.4,2.8);
\draw (7.4,-2.8) -- (7.4,2.8);
%\foreach \x in {(0,0), (7.4,-2.8), (7.4,2.8), (5.1,5.1)}
 %\filldraw \x circle (2pt);
 
\draw[<-,yshift=-0.3cm,line width=0.3pt] (3.92,1.48)--(4.96,1.88); %CE
\draw[->,yshift=-0.3cm,line width=0.3pt] (2.44,0.92)--(3.48,1.32);

\draw[<-,yshift=-0.3cm,line width=0.3pt] (3.81,-1.44)--(4.87,-1.78); %BC
\draw[->,yshift=0.3cm,line width=0.3pt] (3.81,-1.44)--(4.87,-1.78);

\draw[->,xshift=-0.24cm,line width=.3pt] (7.4,-0.49)--(7.4,0.49); %BE
\draw[<-,xshift=0.24cm,line width=.3pt] (7.4,-0.49)--(7.4,0.49);

\end{scope}     
\end{tikzpicture}  \qquad
\begin{tikzpicture}[scale=.25]
\draw (0,0) -- (7.4,-2.8) -- (7.4,2.8) -- (5.1,5.1) -- (0,0);
\draw[dashed] (0,0) -- (7.4,2.8);
\draw (5.1,5.1) -- (7.4,-2.8);
\foreach \x in {(0,0), (7.4,-2.8), (7.4,2.8), (5.1,5.1)}
  \filldraw \x circle (2pt);

\begin{scope}[xshift=7.4cm,yshift=-2.8cm]
\draw (0,0) -- (7.4,-2.8) -- (7.4,2.8) -- (5.1,5.1) -- (0,0);
\draw[dashed] (0,0) -- (7.4,2.8);
\draw (5.1,5.1) -- (7.4,-2.8);
\foreach \x in {(0,0), (7.4,-2.8), (7.4,2.8), (5.1,5.1)}
  \filldraw \x circle (2pt);
\end{scope}

\begin{scope}[xshift=7.4cm,yshift=2.8cm]
\draw (0,0) -- (4.05,-1.55);
\draw[dashed] (4.05,-1.55) -- (7.4,-2.8);
\draw (7.4,-2.8) -- (7.4,2.8) -- (5.1,5.1);
\draw[dashed] (5.1,5.1) -- (0,0) -- (7.4,2.8);
\draw (5.1,5.1) -- (7.4,-2.8);
\foreach \x in {(0,0), (7.4,-2.8), (7.4,2.8), (5.1,5.1)}
  \filldraw \x circle (2pt);
\end{scope}

\begin{scope}[xshift=5.1cm,yshift=5.1cm]
\draw (0,0) -- (7.4,-2.8) -- (7.4,2.8) -- (5.1,5.1) -- (0,0);
\draw[dashed] (0,0) -- (7.4,2.8);
\draw (5.1,5.1) -- (7.4,-2.8);
\foreach \x in {(0,0), (7.4,-2.8), (7.4,2.8), (5.1,5.1)}
  \filldraw \x circle (2pt);
\end{scope}   
\end{tikzpicture} 
\caption{\label{SP8}}
\end{figure}
In general, for any $n \in \mathbb{N}$, we introduce a numbering of the $8\cdot 4^{n-1}$ equilateral triangles of perimeter $2^{1-n}\pi$ which are the boundaries of the faces of all cut out octahedra in the $n$-th step of the construction of the Sierpinski pyramid and denote the numbered triangles 
with $\Delta_{n,m,k}, \, m \in \{ 1,\ldots, 4^{n-1}\}, \, k \in \{1, \ldots , 8\}$.  The numbering of the triangles obeys the following rules: for any $m \in \{1,\ldots, 4^{n-1}\}$ 
\begin{enumerate}
\item $\Delta_{n,m,k}, \, k=1,2,3$ for the boundaries of the holes in an oblique face of $P_{n-1,m},$ 
\item $\Delta_{n,m,4}$ for the boundary of the hole in the horizontal face of $P_{n-1,m}$, and 
\item $\Delta_{n,m,k}, \, k=5,\ldots,8$ for the boundaries of the remaining faces of the octahedron which is cut out of $P_{n-1,m}$.
\end{enumerate}
\end{numbering}
We will construct unbounded Fredholm modules over the algebra $C(SP)$ based on the triangles $\Delta_ {n,m,k}$ similarly as for the Sierpinski gasket.  
We endow each face of any pyramid with the orientation induced by the outer normal vector and orient the boundary correspondingly.   
With this orientation we construct an unbounded Fredholm module  $UFM(\Delta_{0,k})$ over $C(SP)$ in the same way we did for the initial triangle $\Delta_{0,1}$ for the Sierpinski gasket. 
For $n \geq 1$ and $m$ fixed in $\{ 1, \ldots, 4^{n-1}\}$ the triangles  $\Delta_{n,m,k}, \, 1\leq k \leq 8,$ are all boundaries of a face of an octahedron. 
Then for each set of indices $(n,m,k)$ we can construct an unbounded Fredholm module $UFM(\Delta_{n,m,k})$ for $C(SP).$ 
By reversing the orientation in the parameterization of the triangles we obtain the unbounded Fredholm modules denoted by $\overline{UFM(\Delta_{n,m,k})}$. \\
\indent   To obtain a  family of spectral triples which would encode as well as possible the fractal geometry of the pyramid and generate its $K^1$-group we use the same idea as for the Sierpinksi gasket. We first construct a spectral triple   which encodes the fractal geometry of the pyramid and induces the $0$-index map, and thus the $0$-element of $K^1(C(SP))$. To this spectral triple we then add (the appropriate number of) circle spectral triples associated to the boundary of each hole (with in the appropriate orientation) in each oblique face of any small pyramid arising in the construction of the Sierpinski pyramid.\\
\indent A spectral triple which encodes the fractal geometry of the Sierpinski pyramid and induces the $0$-index map is given in the following definition.

\begin{Dfn} \label{SumP} 
The direct sum of unbounded Fredholm modules for the Sierpinski pyramid given by 
\begin{displaymath} 
\underset{k=1}{\overset{4}{\oplus}}UFM(\Delta_{0,k}) \oplus \underset{n=1}{\overset{\infty}{\oplus}} \underset{m=1}{\overset{4^{n-1}}{\oplus}}\underset{k=1}{\overset{8}{\oplus}} UFM(\Delta_{n,m,k})
\end{displaymath} 
is a spectral triple, which is denoted $ZPT$ (the {\em zero pyramid triple}). The Hilbert space and the Dirac operator of this spectral triple are denoted  $H_{ZPT}$ and $D_{ZPT},$ respectively.\end{Dfn}

\begin{Thm}   The bounded Fredholm module coming from the polar
  decomposition of $D_{ZPT}$ induces the trivial element of the
  group $K^1(C(SP))$.  \\
  The spectral triple $ZPT $ has the following
geometric properties:
\begin{itemize}
\item[(i)] The metric induced is the geodesic distance. 
\item[(ii)] The $ZPT$ is summable for any positive $s > 2.$ 
Its zeta-function $\zeta_{ZPT}(s)$ is meromorphic with a simple pole at 2 and is given by $
\zeta_{ZPT}(s) \, =8 \cdot  \frac{(2^s-1)(2^s-2)}{2^s-4}\zeta(s).$ 
\item[(iii)] Let $\mu$ denote the normalized 2-dimensional 
Hausdorff measure on the Sierpinski pyramid. Then for any 
Dixmier trace and any continuous function $g$ in $C(SP)$ we have \begin{displaymath}\Tr_{\omega}\left (|D_{ZPT}|^{-2}\pi_{ZPT}(g)\right) =  \frac{6}{\log 2}\cdot \zeta(2) \cdot \int_{SP}g(x) d\mu(x). \end{displaymath}
\end{itemize}
\end{Thm}
\begin{proof}
Using similar arguments as in \cite{CIL}, pages 27--28, we can  
show that the direct sum is a spectral triple (this result and (i) were also obtained in \cite{Ha}).\\
Let $u$ be a unitary of $C(SP)$. We refer to Proposition \ref{KUFM} (i) to write the K-homology element induced by $ZPT$.
\begin{displaymath}
\Phi_{ZPT}([u])=\underset{k=1}{\overset{4}{\sum}}w_{\Delta_{0,k}}(u) + 
\underset{n=1}{\overset{\infty}{\sum}} \underset{m=1}{\overset{4^{n-1}}{\sum}} \underset{k=1}{\overset{8}{\sum}} w_{\Delta_{n,m,k}} (u). 
\end{displaymath}
Since from a certain number $n_0$ on the winding number of $u$ around any $\Delta_{n,m,k}, \, n>n_0, \, m =1, \ldots, 4^{n-1}, \, k=1, \ldots, 8,$ vanishes, the sum is finite.
Let us compute it.  We notice that $\underset{k=1}{\overset{4}{\sum}}w_{\Delta_{0,k}}(u)=0$ since each edge of the starting tetrahedron is covered twice but in opposite directions and thus the sum of the winding numbers of $u$ around the four triangles $\Delta_{0,k}$ is 0, see Figure 9.  
In fact for any $n \in \{ 1,\ldots, n_0\}$ and $m \in \{1,\ldots 4^{n-1}\}$ the expression $\underset{k=1}{\overset{8}{\sum}} w_{\Delta_{n,m,k}} (u)$ is 0 since each edge in the octahedron numbered $(n,m)$ will be counted twice and in opposite directions. Thus the sum of the winding numbers of $u$ around the eight triangles is 0.  In conclusion the corresponding element in $K^1(C(SP))$ of $ZPT$ is trivial. 
We now fix a positive number $s>2$ and compute the zeta-function of $ZPT$ in $s$.  
We  remark first that, for any $n \in \mathbb{N}_0$, the zeta-function in a point $s>2$ for  $\Delta_{0,k}, \,  k=1,\ldots, 4$, and $\Delta_{n,m,k}, \, n \in \mathbb{N}, \, m=1,\ldots, 4^{n-1}, \, k=1, \dots , 8,$  is the zeta-function in $s$ for $\Delta_{n,j} $ from the Sierpinski gasket case, and thus equals $2 \cdot 2^{-ns} \cdot (2^s-1) \cdot \zeta(s)$, see the equations (\ref{zeta1}) - (\ref{zeta3}).
There are four triangles $\Delta_{0,k}$ and $8\cdot 4^{n-1}$ triangles $\Delta_{n,m,k}$.  Hence the zeta-function for $ZPT$ in a point $s>2$ is
\begin{align*}
\zeta_{ZPT}(s) &= 4\cdot 2(2^s-1)\cdot \zeta(s) +
\sum_{n=1}^{\infty}8\cdot 4^{n-1}\cdot 2\cdot 2^{-ns}(2^s -1) \cdot \zeta(s) \\
&= 8 \cdot  (2^s-1)\zeta(s) +8\cdot 2^{1-s}\cdot (2^s-1)\frac{1}{1 - (4/2^s)} \cdot \zeta(s)\\
&= 8 \cdot  (2^s-1)\zeta(s)\left (1+2^{1-s}\cdot \frac{1}{1-4/2^s}
\right )\\ & =8 \cdot  \frac{(2^s-1)(2^s-2)}{2^s-4}\zeta(s).
\end{align*}
Further on, we obtain 
\begin{align*} 
 \mathrm{Tr}_{\omega}\left (\vert D_{ZPT}\vert ^{- 2}\right  )\,  
  =\, & \underset{x \to 1+}{\lim} (x-1) \zeta_{ZPT}\left ( x\cdot 2 \right )\\ 
= \,& 4\cdot 3 \zeta(2) \underset{x \to 1+}{\lim} (x-1) \frac{1}{4^{x-1}-1}\\
 = \,&   \frac{6}{\log 2} \cdot \zeta(2).
\end{align*}
\end{proof}
In accordance with the convention in Numbering \ref{numb} we shall index the non-horizontal holes of Sierpinski pyramid by 
\begin{displaymath}
HSP:= \{ (n,m,k):\, n\in \mathbb{N}, \, m=1,\ldots, 4^{n-1}, \, k=1,2, 3 \},
\end{displaymath}
and thus we may write
\begin{align*}
K_1(C(SP))& = \oplus_{(n,m,k) \in HSP} \bz,\\
K^1(C(SP)) & =\mathrm{Hom}(K_1(C(SP)), \bz)=\Pi_{(n,m,k) \in HSP} \bz .
\end{align*} 

We shall now construct a spectral triple which induces any prescribed element in  
the group $K^1(C(SP))$.

\begin{Thm}  Let $(f(n,m,k))_{n,m,k} \in \Pi_{(n,m,k) \in HSP}\bz$ then 
\begin{displaymath}
ST(f):=ZPT\oplus
 UFM(f),
 \end{displaymath}
where  $UFM(f)$
  denotes the direct sum \begin{align*}
& \underset{(n,m,k) \in HSP: f(n,m,k) \neq 0}{\oplus} \\
& \begin{cases} UFM(\Delta_{n,m,k})\overset{f(n,m,k)}{\oplus \ldots \oplus} 
UFM(\Delta_{n,m,k}),
  \text{ if } f(n,m,k) > 0. \\  
\overline{UFM(\Delta_{n,m,k})}\overset{-f(n,m,k)}{\oplus \ldots \oplus} 
\overline{UFM(\Delta_{n,m,k})},
  \text{ if } f(n,m,k) < 0, \\  \end{cases}\end{align*}
is a  spectral triple which induces the geodesic distance on the pyramid and the $K^1$-element $(f(n,m,k))_{n,m,k}$. 
\end{Thm}

As for the gasket it is easily seen that if a function $(f(n,m,k))$ is bounded then the summability properties for $ST(f)$ are the same as for $ZPT.$ It is also possible to check that if a function $(f(n,m,k))$ has the property that it is constant for all indices with $n\geq N$ for some natural number $N$ then the volume form will be proportional to that of the 2-dimensional Hausdorff measure.  

At last we will examine the possibility to have an unbounded function 
$f(n,m,k)$ such that the triple $ST(f)$ is summable for any $p>2.$ We
write down - at least formally - the zeta-function $\zeta_f(s)$ for
$ST(f),$ and we obtain:
\begin{align*}
& \zeta_f(s)= \\
& \zeta_{ZPT}(s) + \underset{n=1}{\overset{\infty}{\sum}} \underset{m=1}{\overset{4^{n-1}}{\sum}} \underset{k=1}{\overset{3}{\sum}}\vert f(n,m,k) \vert (2 \cdot 2^{-ns} \cdot (2^s-1) \cdot \zeta(s)) \\
& \zeta_{ZPT}(s)+ 2\cdot (2^s-1) \zeta(s)\underset{n=1}{\overset{\infty}{\sum}} \left( \underset{m=1}{\overset{4^{n-1}}{\sum}} \underset{k=1}
{\overset{3}{\sum}}\vert f(n,m,k) \vert \right )\cdot 2^{-ns}.
\end{align*}

Based on the root criterion we then get:

\begin{Thm}
Let $(f(n,m,k)) \in \Pi_{HSP}\mathbb{Z}$ then $ST(f)$ is summable for any $p>2$ if 
\begin{displaymath}
\underset{n\rightarrow \infty}{\limsup} \left ( \underset{m=1}{\overset{4^{n-1}}{\sum}} \underset{k=1}{\overset{3}{\sum}}\vert f(n,m,k) \vert \right ) ^{1/n} \leq 4.
\end{displaymath}
\end{Thm}

\end{document}